\documentclass{birkart}%
\usepackage{amsmath}
\usepackage[dvips]{graphicx}
\usepackage{amsxtra}
\usepackage{amsfonts}
\usepackage{amssymb}%
\newtheorem{theorem}{Theorem}[section]
\newtheorem{corollary}[theorem]{Corollary}
\newtheorem{lemma}[theorem]{Lemma}

\newtheorem{proposition}[theorem]{Proposition}
\theoremstyle{definition}

\newtheorem{notations}[theorem]{Notations}
\newtheorem{definition}[theorem]{Definition}
\newtheorem{definitions}[theorem]{Definitions}
\newtheorem{remark}[theorem]{Remark}

\newtheorem{observation}[theorem]{Observation}
\newtheorem{example}[theorem]{Example}
\newtheorem{examples}[theorem]{Examples}

\theoremstyle{remark}

\newtheorem*{acknowledgements}{Acknowledgements}
\numberwithin{equation}{section}
\newlength{\displayboxwidth}

\newbox\ipbox
\newcommand{\ip}[2]{\left\langle #1\mathrel{\mathchoice
{\setbox\ipbox=\hbox{$\displaystyle \left\langle\mathstrut #1#2\right\rangle$}
\vrule height\ht\ipbox width0.25pt depth\dp\ipbox}
{\setbox\ipbox=\hbox{$\textstyle \left\langle\mathstrut #1#2\right\rangle$}
\vrule height\ht\ipbox width0.25pt depth\dp\ipbox}
{\setbox\ipbox=\hbox{$\scriptstyle \left\langle\mathstrut #1#2\right\rangle$}
\vrule height\ht\ipbox width0.25pt depth\dp\ipbox}
{\setbox\ipbox=\hbox{$\scriptscriptstyle \left\langle\mathstrut #1#2\right\rangle$}
\vrule height\ht\ipbox width0.25pt depth\dp\ipbox}
} #2\right\rangle}
\def\openone
{\mathchoice
{\hbox{\upshape \small1\kern-3.3pt\normalsize1}}
{\hbox{\upshape \small1\kern-3.3pt\normalsize1}}
{\hbox{\upshape \tiny1\kern-2.3pt\SMALL1}}
{\hbox{\upshape \Tiny1\kern-2pt\tiny1}}}
\renewcommand{\theenumi}{\roman{enumi}}
\renewcommand{\labelenumi}{\textup{(\theenumi )}}
\mathcode`\;="303B
\renewcommand{\subjclassname}{\textup{2000} Mathematics Subject Classification}
\begin{document}
\title[Attractors for representations]{Closed subspaces which are attractors for\\representations of the Cuntz algebras}
\author{Palle E. T. Jorgensen}
\address{Department of Mathematics\\
The University of Iowa\\
14 MacLean Hall\\
Iowa City, IA 52242-1419\\
U.S.A.}
\email{jorgen@math.uiowa.edu}
\thanks{Work supported in part by the U.S. National Science Foundation under grants
DMS-9987777, INT-9724781, and DMS-0139473.}

\begin{abstract}
We analyze the structure of co-invariant subspaces for representations of the
Cuntz algebras $\mathcal{O}_{N}$ for $N=2,3,\dots$, $N<\infty$, with special
attention to the representations which are associated to orthonormal and
tight-frame wavelets in $L^{2}\left(  \mathbb{R}\right)  $ corresponding to
scale number $N$.

\end{abstract}
\maketitle
\tableofcontents

\section{\label{Int}Introduction: Wavelets}

\setlength{\displayboxwidth}{\textwidth}\addtolength{\displayboxwidth
}{-2\leftmargini}A particular construction of wavelets on the real line
$\mathbb{R}$ is based on what is called subband filters. The idea is that a
wavelet decomposition of $L^{2}\left(  \mathbb{R}\right)  $ can be organized
in frequency bands with adjustment to a system of subspaces (identifying a
cascade of resolutions) in $L^{2}\left(  \mathbb{R}\right)  $, and each
frequency band having its scaling resolution. The subband filters may be
realized as functions on the torus $\mathbb{T}=\left\{  z\in\mathbb{C}%
:\left\vert z\right\vert =1\right\}  $. If the scaling is $N>1$ then the
$L^{2}\left(  \mathbb{R}\right)  $ functions needed in the wavelet
decomposition may be obtained (under favorable conditions) as solutions to a
system of cocycle conditions, see (\ref{eqInt.1}) below, involving Fourier
transform $\psi\mapsto\hat{\psi}$ on $\mathbb{R}$, or rather $L^{2}\left(
\mathbb{R}\right)  $. The $L^{2}\left(  \mathbb{R}\right)  $-system consists
of a scaling function $\varphi$ and the wavelet generators $\psi_{1}%
,\dots,\psi_{N-1}$. The cocycle conditions (see \cite{BrJo02b} for details)
are
\begin{equation}
\begin{aligned} \sqrt{N}\hat{\varphi}\left( N\xi\right) & =m_{0}\left( e^{i\xi}\right) \hat{\varphi}\left( \xi\right) ,\qquad\xi\in\mathbb{R}\text{, and} \\ \sqrt{N}\hat{\psi}_{j}\left( N\xi\right) & =m_{j}\left( e^{i\xi}\right) \hat{\varphi}\left( \xi\right) ,\qquad j=1,\dots,N-1. \end{aligned} \label{eqInt.1}%
\end{equation}
For more details on this, we refer the reader to
\cite{Dau92}, \cite{Mal99},
\cite{Jor99b}, \cite{BrJo99b}%
, \cite{Jor01b}, and \cite{BrJo02b}. After the problem is discretized, and a
Fourier series is introduced, we then arrive at a certain system of operators
on $L^{2}\left(  \mathbb{T}\right)  $, where the one-torus $\mathbb{T}$ is
equipped with the usual normalized Haar measure. The operators are defined
from a fixed system of functions $m_{0},\dots,m_{N-1}$ on $\mathbb{T}$ as
follows:%
\begin{equation}
\left(  S_{j}f\right)  \left(  z\right)  =m_{j}\left(  z\right)  f\left(
z^{N}\right)  ,\qquad f\in L^{2}\left(  \mathbb{T}\right)  ,\;z\in
\mathbb{T},\;j=0,1,\dots,N-1. \label{eqInt.2}%
\end{equation}
The orthogonality conditions which are usually imposed are known to imply the
following relations on the operators $\left(  S_{j}\right)  _{j=0}^{N-1}$ in
(\ref{eqInt.2}):%
\begin{equation}
S_{j}^{\ast}S_{i}=\delta_{j,i}I_{L^{2}\left(  \mathbb{T}\right)  }\text{\quad
and\quad}\sum_{j=0}^{N-1}S_{j}S_{j}^{\ast}=I_{L^{2}\left(  \mathbb{T}\right)
}. \label{eqInt.3}%
\end{equation}
The reader is referred to \cite{Jor01b}, \cite{Jor00b}, and \cite{BrJo02b} for
additional details on this point. Specifically, the properties which must be
imposed on the functions $m_{0},\dots,m_{N-1}$ from (\ref{eqInt.2}) are known
in signal processing as the subband quadrature (if $N=2$) conditions. They can
be checked to be equivalent to the operator relations (\ref{eqInt.3}), also
known as the Cuntz relations. They are satisfied if and only if the $N\times
N$ matrix ($i,j=0,\dots,N-1$)%
\begin{equation}
A_{i,j}\left(  z\right)  :=\frac{1}{N}\sum_{\substack{w\in\mathbb{T}\\w^{N}%
=z}}m_{i}\left(  w\right)  w^{-j},\qquad z\in\mathbb{T}, \label{eqInt.4}%
\end{equation}
is unitary for all (or almost all, with respect to Haar measure on
$\mathbb{T}$) $z\in\mathbb{T}$. The relations (\ref{eqInt.3}) are called the
Cuntz relations, and they are special cases of relations which are defined
axiomatically in the theory of representation of $C^{\ast}$-algebras. But they
have an independent life in the science of signal processing; see
\cite{BrJo02b}, \cite{Mal99}, and the references given there.

\section{\label{Cun}The Cuntz relations}

It is known \cite{Cun77} that there is a simple $C^{\ast}$-algebra
$\mathcal{O}_{N}$ such that the representations $\rho$ of $\mathcal{O}_{N}$
are in a $1$--$1$ correspondence with systems of operators on Hilbert space
satisfying the Cuntz relations. If the generators of $\mathcal{O}_{N}$ are
denoted $s_{i}$, then
\begin{equation}
s_{i}^{\ast}s_{j}=\delta_{i,j}1\text{\quad and\quad}\sum_{j=1}^{N}s_{j}%
s_{j}^{\ast}=1, \label{eqCun.1}%
\end{equation}
where $1$ denotes the unit-element in the $C^{\ast}$-algebra $\mathcal{O}_{N}%
$. The system of operators $S_{j}:=\rho\left(  s_{j}\right)  $ will then
satisfy the Cuntz relations on the Hilbert space $\mathcal{H}$ which carries
the representation $\rho$. Conversely, every system of operators $S_{j}$ on a
Hilbert space $\mathcal{H}$ which satisfies the Cuntz relations comes from a
representation $\rho$ of $\mathcal{O}_{N}$ via the formula $\rho\left(
s_{j}\right)  =S_{j}$, $j=1,\dots,N$. It is known that not every
representation of $\mathcal{O}_{N}$ is of the form (\ref{eqInt.2}).
$\mathcal{O}_{N}$ has many type $\operatorname*{III}$ representations, and
(\ref{eqInt.2}) implies type $\operatorname*{I}$ (see \cite{BJO99}). But it
was shown in \cite{Jor01b} that the analysis of the wavelet representations
(\ref{eqInt.2}) predicts a number of global and geometric properties of the
variety of all wavelets subject to a fixed scaling.

In addition the representations of $\mathcal{O}_{N}$ and their functional
models are used in multivariable scattering theory: see, for example,
\cite{DKS01}, \cite{BaVi02a}, \cite{BaVi02b}, and \cite{Kribs}.

While our present results apply to general representations of $\mathcal{O}%
_{N}$, they are motivated by (\ref{eqInt.2}). In particular, we will be
interested in smaller subspaces of $L^{2}\left(  \mathbb{T}\right)  $ which
determine the representation.
Further, when a wavelet representation is
given, a result in Section \ref{Fin} makes
precise a sense in which these
subspaces are attractors, i.e., they
arise as limits of a dynamical iteration.
It is known from \cite{Jor01b} that, if the
wavelet system consists of compactly supported functions $\varphi,\psi
_{1},\dots,\psi_{N-1}$ on $\mathbb{R}$, then there is a finite-dimensional
subspace $\mathcal{L}$ which determines the representation, and therefore the
wavelet analysis. The construction of wavelets from multiresolutions (i.e.,
scales of closed subspaces in $L^{2}\left(  \mathbb{R}\right)  $) is explained
in several books on wavelets, starting with \cite{Dau92}. A geometric approach
which is close to the present one, based on wandering subspaces, was first
outlined in \cite{DaLa98}; see also \cite{BaMe99}. In this paper, we study
such determining minimal closed subspaces in the context of the most general
representations of the Cuntz algebras $\mathcal{O}_{N}$.

\section{\label{Sub}Subspaces of the Hilbert space}

The Hilbert space which carries a representation of $\mathcal{O}_{N}$ must be
in\-fi\-nite-di\-men\-sion\-al, but we show that it contains \textquotedblleft
small\textquotedblright\ distinguished closed subspaces.

First some definitions: there will be eleven in all.

\begin{definitions}
\label{DefSub.1}

\begin{enumerate}
\item \label{DefSub.1(0)}If $\left\{  S_{i}\right\}  _{i=1}^{N}$ is a
representation of $\mathcal{O}_{N}$ on a Hilbert space, set $S_{I}=S_{i_{1}%
}\cdots S_{i_{k}}$ for all \emph{multi-indices} $I=\left(  i_{1},\dots
,i_{k}\right)  $, and define the \emph{multi-index length} as $\left\vert
I\right\vert :=k$. The set of all such multi-indices $I$ with $i_{\nu}%
\in\left\{  1,\dots,N\right\}  $ is denoted $\mathcal{I}\left(  N\right)  $.
The set of multi-indices of length $k$ is denoted $\mathcal{I}_{k}\left(
N\right)  $.

\item \label{DefSub.1(1)}If $\mathcal{F}$ is a set of vectors in a Hilbert
space $\mathcal{H}$, the notation $\bigvee\mathcal{F}$ stands for the closed
linear span of $\mathcal{F}$.

\item \label{DefSub.1(2)}If $\mathcal{M}$ is a family of subspaces in
$\mathcal{H}$, then $\bigvee\mathcal{M}$ stands for the closed linear span of
these subspaces.

\item \label{DefSub.1(3)}If $\left\{  S_{i}\right\}  $ are the isometries
which define a representation of $\mathcal{O}_{N}$, acting on a fixed Hilbert
space $\mathcal{H}$, and if $\mathcal{L}$ is a closed subspace in
$\mathcal{H}$, then we set
\[
S\mathcal{L}:=\bigvee_{i=1}^{N}S_{i}\mathcal{L}\text{,\quad and\quad}S^{\ast
}\mathcal{L}:=\bigvee_{i=1}^{N}S_{i}^{\ast}\mathcal{L}.
\]
(Note that, since the operators $S_{i}$ are isometries, it follows that each
space $S_{i}\mathcal{L}$ is closed. From the identity%
\[
\sum_{i=1}^{N}S_{i}S_{i}^{\ast}=I_{\mathcal{H}},
\]
it follows that there are natural conditions which imply that each of the
spaces $S_{i}^{\ast}\mathcal{L}$ is also closed. We spell out this point in
Lemma \textup{\ref{LemSub.2}} below.)

\item \label{DefSub.1(4)}A closed subspace $\mathcal{L}$ is said to be
\emph{co-invariant} for a fixed representation $\left(  S_{i}\right)  $ of
$\mathcal{O}_{N}$ acting in a Hilbert space $\mathcal{H}$ if $S^{\ast
}\mathcal{L}\subset\mathcal{L}$.

\item \label{DefSub.1(5)}A closed subspace $\mathcal{L}$ is said to be
\emph{saturated} if
\[
\bigvee_{k=1}^{\infty}S^{k}\mathcal{L}=\mathcal{H},
\]
where, for every $k$,%
\[
S^{k}\mathcal{L}:=\!\!\bigvee_{\left(  i_{1},\dots,i_{k}\right)
\in\mathcal{I}_{k}\left(  N\right)  }\!\!S_{i_{1}}S_{i_{2}}\cdots S_{i_{k}%
}\mathcal{L}.
\]

\item \label{DefSub.1(New6)}A co-invariant closed subspace $\mathcal{L}$ is
said to be \emph{minimal} if the corresponding complementing space
$\mathcal{W}:=S\mathcal{L}\ominus\mathcal{L}$ generates a maximal subspace,
i.e., the family
\[
\mathcal{F}\left(  \mathcal{W}\right)  :=\left\{  \mathcal{W},\;S_{I}%
\mathcal{W}:I=\left(  i_{1},\dots,i_{k}\right)  ,\;k\geq1,\;I\in
\mathcal{I}\left(  N\right)  \right\}
\]
is such that%
\[
\sideset{}{^{\oplus}}{\displaystyle\sum}\left\{  \mathcal{K}%
:\mathcal{K}\in\mathcal{F}\left(  \mathcal{W}\right)  \right\}
\]
is maximal in $\mathcal{H}$.

\item \label{DefSub.1(6)}We say that $\mathcal{L}\subset\mathcal{H}$ is a
\emph{core} for the representation if it is co-invariant, saturated, and minimal.

\item \label{DefSub.1(6bis)}A closed subspace $\mathcal{W}\subset\mathcal{H}$
is called \emph{wandering} if all the spaces in the family $\mathcal{F}\left(
\mathcal{W}\right)  $ are mutually orthogonal.

\item \label{DefSub.1(7)}A subspace \textup{(}not necessarily
closed\/\textup{)} $\mathcal{V}$ in $\mathcal{H}$ is said to \emph{reduce to a
closed subspace} $\mathcal{L}$ \emph{in a representation} $\left(
S_{i}\right)  $ if, for every $v\in\mathcal{V}$, there is a $k_{0}%
\in\mathbb{N}$ such that $S_{i_{1}}^{\ast}S_{i_{2}}^{\ast}\cdots S_{i_{k}%
}^{\ast}v\in\mathcal{L}$ whenever $k\geq k_{0}$.

\item \label{DefSub.1(8)}A subspace $\mathcal{L}$ is said to be \emph{stable}
with respect to a representation $\left(  S_{i}\right)  $ of $\mathcal{O}_{N}$
if it is invariant under each of the projections $E_{i}:=S_{i}S_{i}^{\ast}$,
$i=1,\dots,N$.
\end{enumerate}
\end{definitions}

Our main result is Theorem \ref{ThmFin.1}.
It states that the standard wavelet
representation (\ref{eqInt.2}) corresponding to
filter functions $m_{j}$ which are Lipschitz
always has a natural finite-dimensional
co-invariant subspace $\mathcal{L}_{\operatorname*{fin}}$. The
results leading up to Section \ref{Fin} throw
light on the properties of the
subspace $\mathcal{L}_{\operatorname*{fin}}$, and on co-invariant
subspaces more generally.

\begin{lemma}
\label{LemSub.2}Let a representation $\left(  S_{i}\right)  $ of
$\mathcal{O}_{N}$ on a Hilbert space $\mathcal{H}$ be given. Let
$\mathcal{L}\subset\mathcal{H}$ be a closed subspace which is stable. Then
each of the linear spaces $S_{i}^{\ast}\mathcal{L}=\left\{  S_{i}^{\ast}%
x:x\in\mathcal{L}\right\}  $ is closed in $\mathcal{H}$.
\end{lemma}

\begin{proof}
Set $i=1$ for specificity. Let $x_{n}$ be a sequence in $\mathcal{L}$ such
that $S_{1}^{\ast}x_{n}\underset{n\rightarrow\infty}{\longrightarrow}y$. We
introduce the projection $E_{1}=S_{1}S_{1}^{\ast}$. Since $S_{j}^{\ast}%
E_{1}=\delta_{j,1}S_{1}^{\ast}$, we have $S_{1}^{\ast}E_{1}x_{n}=S_{1}^{\ast
}x_{n}$, $S_{j}^{\ast}E_{1}x_{n}=0$ if $j\neq1$, and $\left\Vert S_{1}^{\ast
}E_{1}\left(  x_{n}-x_{m}\right)  \right\Vert =\left\Vert E_{1}\left(
x_{n}-x_{m}\right)  \right\Vert \underset{n,m\rightarrow\infty}%
{\longrightarrow}0$. The sequence $\left(  E_{1}x_{n}\right)  $ is convergent,
and its limit is in $\mathcal{L}$ since $\mathcal{L}$ is closed and stable. If
$E_{1}x_{n}\underset{n\rightarrow\infty}{\longrightarrow}z$, then $S_{1}%
^{\ast}z=y$, and we conclude that $S_{1}^{\ast}\mathcal{L}$ is closed.
\end{proof}

\section{\label{Co}Co-invariant closed subspaces}

A representation of $\mathcal{O}_{N}$ on a Hilbert space $\mathcal{H}$ is
specified by operators $S_{1},\dots,S_{N}$ on $\mathcal{H}$ subject to the
Cuntz relations
\begin{equation}
S_{i}^{\ast}S_{j}=\delta_{i,j}I_{\mathcal{H}},\qquad\sum_{j=1}^{N}S_{j}%
S_{j}^{\ast}=I_{\mathcal{H}}, \label{eqCo.1}%
\end{equation}
where $I_{\mathcal{H}}$ denotes the identity operator in the Hilbert space
$\mathcal{H}$. Intrinsic to this is the set of $N$ commuting projections
$E_{j}:=S_{j}S_{j}^{\ast}$. They enter into the statement of the next lemmas.

In the next lemma we record some general properties about co-invariant
subspaces. They are stated in terms of projections. Recall there is a $1$--$1$
correspondence between closed subspaces $\mathcal{L}\subset\mathcal{H}$ and
\emph{projections} $P$ in $\mathcal{H}$, i.e., $P=P^{\ast}=P^{2}$. If
$\mathcal{L}$ is given, there is a unique $P$ such that $\mathcal{L}%
=P\mathcal{H}=\left\{  x\in\mathcal{H}:Px=x\right\}  $, and conversely.

\begin{lemma}
\label{LemCoNew.1}Let $\left\{  S_{i}\right\}  _{i=1}^{N}$ be a representation
of $\mathcal{O}_{N}$ on a Hilbert space $\mathcal{H}$, and let $\mathcal{L}$
be a closed subspace in $\mathcal{H}$. Set $\alpha\left(  A\right)
=\sum_{i=1}^{N}S_{i}AS_{i}^{\ast}$, $A\in\mathcal{B}\left(  \mathcal{H}\right)
$.\renewcommand{\theenumi}{\alph{enumi}}

\begin{enumerate}
\item \label{LemCoNew.1(1)}Then $\alpha\colon\mathcal{B}\left(  \mathcal{H}%
\right)  \rightarrow\mathcal{B}\left(  \mathcal{H}\right)  $ is an
endomorphism satisfying $\alpha\left(  I_{\mathcal{H}}\right)  =I_{\mathcal{H}%
}$.

\item \label{LemCoNew.1(2)}If $P$ denotes the projection onto $\mathcal{L}$,
then $\mathcal{L}$ is co-invariant if and only if $P\leq\alpha\left(
P\right)  $; and $\alpha\left(  P\right)  $ is then the projection onto
$S\mathcal{L}$.

\item \label{LemCoNew.1(3)}If $\mathcal{L}$ is co-invariant, then
$Q=\alpha\left(  P\right)  -P$ is a projection. Its range is the wandering
subspace $\mathcal{W}:=\left(  S\mathcal{L}\right)  \ominus\mathcal{L}$.
\end{enumerate}
\end{lemma}

\begin{proof}
The details are left to the reader. They are based on standard geometric facts
about projections in Hilbert space.
\end{proof}

\begin{lemma}
\label{LemCoNew.2}Let $\left\{  S_{i}\right\}  _{i=1}^{N}$ be a representation
of $\mathcal{O}_{N}$, and let $\mathcal{L}$ be a closed co-invariant subspace.
Set $\mathcal{W}=\left(  S\mathcal{L}\right)  \ominus\mathcal{L}$, and
$\mathcal{F}\left(  \mathcal{W}\right)  =\mathcal{W}\oplus S\mathcal{W}\oplus
S^{2}\mathcal{W}\oplus\cdots$. Let $P_{\mathcal{L}}$ be the projection onto
$\mathcal{L}$, and let $\alpha$ be the endomorphism of Lemma
\textup{\ref{LemCoNew.1}}. Then the limit $P_{\alpha}=\lim_{n\rightarrow
\infty}\alpha^{n}\left(  P_{\mathcal{L}}\right)  $ exists, and the projection
$P_{\mathcal{F}}$ onto $\mathcal{F}\left(  \mathcal{W}\right)  $ is given by
the formula%
\begin{equation}
P_{\mathcal{F}}=P_{\infty}-P_{\mathcal{L}}. \label{eqCoNew.star}%
\end{equation}
Moreover, the operators $T_{i}=P_{\mathcal{F}}S_{i}P_{\mathcal{F}}$,
$i=1,\dots,N$, satisfy the following \textup{(}%
Cuntz--Toeplitz--Fock\/\textup{)} relations:\renewcommand{\theenumi}{\alph{enumi}}

\begin{enumerate}
\item \label{LemCoNew.2(1)}$T_{i}^{\ast}T_{j}=\delta_{i,j}P_{\mathcal{F}}$,
$i,j=1,\dots,N$,

\item \label{LemCoNew.2(2)}$\sum_{i=1}^{N}T_{i}T_{i}^{\ast}\leq P_{\mathcal{F}%
}$,

\item \label{LemCoNew.2(3)}$T_{i}^{\ast}w=0$ for all $i=1,\dots,N$, and all
$w\in\mathcal{W}$.
\end{enumerate}
\end{lemma}

\begin{proof}
The details amount to direct verifications and are left for the reader.
\end{proof}

\begin{lemma}
\label{LemCoNew.3}Let $\left\{  S_{i}\right\}  _{i=1}^{N}$ be a representation
of $\mathcal{O}_{N}$ on a Hilbert space $\mathcal{H}$. Let $\mathcal{W}%
\subset\mathcal{H}$ be a closed subspace such that $\left\langle x\mid
S_{J}y\right\rangle =0$ for all $x,y\in\mathcal{W}$ and all multi-indices
$J\in\mathcal{I}\left(  N\right)  $. Then the operators $T_{i}:=P_{\mathcal{F}%
\left(  \mathcal{W}\right)  }S_{i}P_{\mathcal{F}\left(  \mathcal{W}\right)  }$
satisfy the conditions \textup{(\ref{LemCoNew.2(1)})--(\ref{LemCoNew.2(3)})}
in Lemma \textup{\ref{LemCoNew.2}}, i.e., the subspace $\mathcal{F}\left(
\mathcal{W}\right)  =\mathcal{W}\oplus S\mathcal{W}\oplus S^{2}\mathcal{W}%
\oplus\cdots$ induces a Fock-space representation of the Cuntz--Toeplitz relations.
\end{lemma}

\begin{proof}
It suffices to show that each operator $S_{i}^{\ast}$ maps $\mathcal{W}$ into
$\mathcal{L}=\mathcal{H}\ominus\mathcal{F}\left(  \mathcal{W}\right)  $.
Suppose $x\in\mathcal{W}$ and $y\in\mathcal{F}\left(  \mathcal{W}\right)  $:
then we claim that $\left\langle S_{i}^{\ast}x\mid y\right\rangle =0$. The
assertion follows from this. To establish the claim, suppose first that
$y\in\mathcal{W}$. Then $\left\langle S_{i}^{\ast}x\mid y\right\rangle
=\left\langle x,S_{i}y\right\rangle =0$ holds on account of the definition of
$\mathcal{W}$. Similarly $\left\langle S_{i}^{\ast}x\mid S_{J}y\right\rangle
=\left\langle x\mid S_{i}S_{J}y\right\rangle =0$ for all $x,y\in\mathcal{W}$,
and all multi-indices $J$. It follows, in particular, that the operators
$T_{i}=P_{\mathcal{F}\left(  \mathcal{W}\right)  }S_{i}P_{\mathcal{F}\left(
\mathcal{W}\right)  }$ satisfy $T_{i}^{\ast}x=0$ for all $x\in\mathcal{W}$, in
addition to the Cuntz--Toeplitz relations (\ref{LemCoNew.2(1)}%
)--(\ref{LemCoNew.2(3)}).
\end{proof}

Our next result is a partial converse to Lemma \ref{LemCoNew.2}.

\begin{proposition}
\label{ThmCoNew.4}Let $\left\{  S_{i}\right\}  _{i=1}^{N}$ be a representation
of $\mathcal{O}_{N}$ on a Hilbert space $\mathcal{H}$, and let $\mathcal{W}%
_{m}$ be chosen as in Lemma \textup{\ref{LemCoNew.6}} below to be maximal
\textup{(}in the sense of Zorn's lemma\/\textup{)} with respect to the
property%
\begin{equation}
\left\langle x\mid S_{J}y\right\rangle =0,\qquad x,y\in\mathcal{W}_{m}\text{,
and }J\in\mathcal{I}\left(  N\right)  . \label{eqCoNew.A1}%
\end{equation}
Then every co-invariant closed subspace $\mathcal{L}$ such that
$S^{\ast}\mathcal{W}_{m}\subset\mathcal{L}\subset
\mathcal{W}_{m}^{\perp}$ satisfies%
\begin{equation}
S\mathcal{L}\ominus\mathcal{L}=\mathcal{W}_{m}. \label{eqCoNew.A2}%
\end{equation}
Conversely, if equality holds in \textup{(\ref{eqCoNew.A2})} for some
wandering subspace $\mathcal{W}$ and all $\mathcal{L}$ with
$S^{\ast}\mathcal{W}\subset\mathcal{L}\subset
\mathcal{W}^{\perp}$, then
$\mathcal{W}$ is maximal.
\end{proposition}

\begin{proof}
First note that, if $\mathcal{W}$ is any wandering
subspace, then the inclusions
\[
S^{\ast}\mathcal{W}\subset\mathcal{F}\left( \mathcal{W}\right) ^{\perp}\subset
\mathcal{W}^{\perp}
\]
are automatic; and
$\mathcal{L}_\mathcal{W}:=\mathcal{F}\left( \mathcal{W}\right) ^{\perp}$ is
co-invariant. Also $\mathcal{L}_\mathcal{W}$ is maximal
among the co-invariant subspaces $\mathcal{L}$ satisfying
\[
S^{\ast}\mathcal{W}\subset\mathcal{L}\subset
\mathcal{W}^{\perp}.
\]

Recall that, if $P_{\mathcal{L}}$ is the projection onto $\mathcal{L}$, then
by Lemma \ref{LemCoNew.1}, $\alpha\left(  P_{\mathcal{L}}\right)  =\sum
_{i=1}^{N}S_{i}P_{\mathcal{L}}S_{i}^{\ast}$ is the projection onto
$S\mathcal{L}$. If we show that $\alpha\left(  P_{\mathcal{L}}\right)  x=x$
for all $x\in\mathcal{W}_{m}$, it follows that $\mathcal{W}_{m}\subseteq
\left(  S\mathcal{L}\right)  \ominus\mathcal{L}$. But the vectors in $\left(
S\mathcal{L}\right)  \ominus\mathcal{L}$ satisfy the orthogonality relations
(\ref{eqCoNew.A1}), and we may then conclude that (\ref{eqCoNew.A2}) holds by
virtue of the maximality of $\mathcal{W}_{m}$; see Lemma \ref{LemCoNew.6}
below. The proof is now completed since $S_{i}^{\ast}x\in\mathcal{L}$ for all
$x\in\mathcal{W}_{m}$, and we get $\alpha\left(
P_{\mathcal{L}}\right)  x=\sum_{i=1}^{N}S_{i}P_{\mathcal{L}}S_{i}^{\ast}%
x=\sum_{i=1}^{N}S_{i}S_{i}^{\ast}x=x$. We leave the verification of the
converse implication to the reader, i.e., that equality in (\ref{eqCoNew.A2})
for some $\mathcal{W}$ and all $\mathcal{L}$ as specified implies that
$\mathcal{W}$ is maximal with respect to (\ref{eqCoNew.A1}).
\end{proof}

\begin{remark}
\label{RemCo.minushalf}The argument shows in particular that if
$\mathcal{W}$ is any closed wandering subspace, then
$S^{\ast}\mathcal{W}\subset
\mathcal{W}^{\perp}$, and further that there will always be
co-invariant subspaces $\mathcal{L}$ satisfying
\begin{equation}
S^{\ast}\mathcal{W}\subset\mathcal{L}\subset
\mathcal{W}^{\perp}.
\label{eqCoJan29.pound}
\end{equation}
In fact, if $\mathcal{W}$ is given, then
$\mathcal{L}=\mathcal{F}\left( \mathcal{W}\right) ^{\perp}$ is
such an intermediate co-invariant subspace. If $\mathcal{W}$
is also assumed maximal in the sense of Proposition \textup{\ref{ThmCoNew.4}},
then it follows that
$\mathcal{L}=\mathcal{F}\left( \mathcal{W}\right) ^{\perp}$
is the only closed and
co-invariant saturated subspace which is intermediate as
stated in \textup{(\ref{eqCoJan29.pound})}.
\end{remark}

\begin{remark}
\label{RemCo.0}While, on the face of it, the relations \textup{(\ref{eqCo.1})}
appear to depend on a choice of coordinates, this is not so: the Cuntz
$C^{\ast}$-algebra is in fact a functor from Hilbert space to $C^{\ast}%
$-algebras, as noted in \cite{Cun77}. Specifically, consider the mapping
\begin{equation}
z=\left(  z_{1},\dots,z_{N}\right)  \longmapsto s\left(  z\right)  =\sum
_{j=1}^{N}z_{j}s_{j} \label{eqRemCo.0}%
\end{equation}
from the Hilbert space $\mathbb{C}^{N}$ into the $C^{\ast}$-algebra
$\mathcal{O}_{N}$. If $1$ denotes the unit element in $\mathcal{O}_{N}$, we
get $s\left(  z\right)  ^{\ast}s\left(  w\right)  =\left\langle z\mid
w\right\rangle 1$ where $\left\langle z\mid w\right\rangle =\sum_{j=1}^{N}%
\bar{z}_{j}w_{j}$ is the usual inner product of $\mathbb{C}^{N}$. If
$\left\Vert \,\cdot\,\right\Vert $ denotes the $C^{\ast}$-norm on
$\mathcal{O}_{N}$, i.e., satisfying $\left\Vert 1\right\Vert =1$ and
$\left\Vert A^{\ast}A\right\Vert =\left\Vert A\right\Vert ^{2}$,
$A\in\mathcal{O}_{N}$, then $\left\Vert s\left(  z\right)  \right\Vert
=\left(  \sum_{j=1}^{N}\left\vert z_{j}\right\vert ^{2}\right)  ^{\frac{1}{2}%
}=\left\Vert z\right\Vert _{\mathbb{C}^{N}}$.

\emph{States} on $\mathcal{O}_{N}$ are linear functionals $\omega
\colon\mathcal{O}_{N}\rightarrow\mathbb{C}$ such that $\omega\left(  1\right)
=1$, $\omega\left(  A^{\ast}A\right)  \geq0$, $A\in\mathcal{O}_{N}$. A special
state $\omega$ on $\mathcal{O}_{N}$ is determined by%
\begin{equation}
\omega\left(  s_{i_{1}}s_{i_{2}}\cdots s_{i_{k}}s_{j_{l}}^{\ast}\cdots
s_{j_{2}}^{\ast}s_{j_{1}}^{\ast}\right)  =\prod_{s}\delta_{i_{s},j_{s}}%
\cdot\delta_{k,l}\cdot N^{-k}. \label{eqRemCo.0starstar}%
\end{equation}

As is well known, states induce cyclic representations. A representation
$\rho$ of $\mathcal{O}_{N}$ in a Hilbert space $\mathcal{H}$ is \emph{cyclic}
if there is a vector $\Omega$ in $\mathcal{H}$ such that $\left\{  \rho\left(
A\right)  \Omega:A\in\mathcal{O}_{N}\right\}  $ is dense in $\mathcal{H}$. If
the state is $\omega$, the corresponding representation $\rho$ is determined
by $\omega\left(  A\right)  =\left\langle \Omega\mid\rho\left(  A\right)
\Omega\right\rangle $. When $\rho$ is given, we set $S_{i}:=\rho\left(
s_{i}\right)  $, $i=1,\dots,N$. We shall often identify a representation
$\rho$ of $\mathcal{O}_{N}$ with the corresponding operators $S_{i}%
=\rho\left(  s_{i}\right)  $.
\end{remark}

The representation $\left\{  S_{i}\right\}  _{i=1}^{N}$ of the state
(\ref{eqRemCo.0starstar}) does not have a closed co-invariant and saturated
subspace which is minimal with respect to the usual ordering of closed
subspaces. Nonetheless the following general result holds.

\begin{lemma}
\label{LemCoNew.6}Let $\left\{  S_{i}\right\}  _{i=1}^{N}$ be a representation
of $\mathcal{O}_{N}$ on a Hilbert space $\mathcal{H}$.

\begin{enumerate}
\item \label{LemCoNew.6(1)}Then $\mathcal{H}$ contains a wandering subspace
$\mathcal{W}_{m}$ which is maximal with respect to inclusion. Specifically,
$\mathcal{W}_{m}$ satisfies the implication: Whenever $\mathcal{W}$ is a
wandering subspace such that $\mathcal{W}_{m}\subseteq\mathcal{W}$, then
$\mathcal{W}_{m}=\mathcal{W}$.

\item \label{LemCoNew.6(2)}Maximal wandering subspaces are non-unique.
\end{enumerate}
\end{lemma}

\begin{proof}
Recall that a closed subspace $\mathcal{W}$ is wandering if and only if the
family of spaces $\left\{  \mathcal{W},\;S_{I}\mathcal{W}\right\}  $ is
orthogonal as $I$ varies over all the multi-indices, i.e., $I\in
\mathcal{I}\left(  N\right)  $. Equivalently,
\begin{equation}
P_{\mathcal{W}}S_{I}P_{\mathcal{W}}=0\text{\qquad for all }I\in\mathcal{I}%
\left(  N\right)  \label{eqLemCoNew.6starstarstar}%
\end{equation}
where $P_{\mathcal{W}}$ denotes the projection onto $\mathcal{W}$. If
$\left\{  \mathcal{W}_{\alpha}\right\}  $ is a linearly ordered family of
wandering subspaces (i.e., for all $\alpha$, $\beta$, one of the inclusions
$\mathcal{W}_{\alpha}\subseteq\mathcal{W}_{\beta}$ or $\mathcal{W}_{\beta
}\subseteq\mathcal{W}_{\alpha}$ holds), then it can be checked that
$\mathcal{W}=\bigvee_{\alpha}\mathcal{W}_{\alpha}$ satisfies the condition
(\ref{eqLemCoNew.6starstarstar}). The existence of $\mathcal{W}_{m}$ now
follows from Zorn's lemma. The example of (\ref{eqRemCo.0starstar}) shows that
$\mathcal{W}_{m}$ will typically be non-unique.
\end{proof}

\begin{lemma}
\label{LemCo.0star}Let $\left(  S_{i}\right)  _{i=1}^{N}$ be a representation
of $\mathcal{O}_{N}$ on a Hilbert space $\mathcal{H}$, and let $\mathcal{L}$
be a closed subspace of $\mathcal{H}$. Then conditions
\textup{(\ref{LemCo.0star(1)})} and \textup{(\ref{LemCo.0star(2)})} are equivalent:

\begin{enumerate}
\item \label{LemCo.0star(1)}$S^{\ast}\mathcal{L}\subseteq\mathcal{L}$;

\item \label{LemCo.0star(2)}$\mathcal{L}\subseteq S\mathcal{L}$.
\end{enumerate}
\end{lemma}

\begin{proof}
(\ref{LemCo.0star(1)})${}\Rightarrow{}$(\ref{LemCo.0star(2)}). Assume
(\ref{LemCo.0star(1)}). Let $x\in\mathcal{L}$. Then $x=\sum_{i}S_{i}%
S_{i}^{\ast}x\in S\mathcal{L}$ since $S_{i}^{\ast}x\in\mathcal{L}$ for
$i=1,\dots,N$.

(\ref{LemCo.0star(2)})${}\Rightarrow{}$(\ref{LemCo.0star(1)}). Assume
(\ref{LemCo.0star(2)}). If $x\in\mathcal{L}$, then there are $y_{i}%
\in\mathcal{L}$ such that $x=\sum_{i}S_{i}y_{i}$. We get $S_{j}^{\ast}%
x=\sum_{i}S_{j}^{\ast}S_{i}y_{i}=\sum_{i}\delta_{j,i}y_{i}=y_{j}$, which
proves (\ref{LemCo.0star(1)}).
\end{proof}

\begin{remark}
\label{RemCo.query}Example \textup{\ref{ExaRes.2}} below shows that there are
representations $\left(  S_{i}\right)  _{i=1}^{N}$ of $\mathcal{O}_{N}$ for
every $N=2,3,\dots$, such that the second inclusion in the lemma is sharp, but
not the first; i.e., it may happen that $S^{\ast}\mathcal{L}=\mathcal{L}$
while $\mathcal{L}\subsetneqq S\mathcal{L}$.
\end{remark}

\begin{example}
\label{ExaCo.0starstarstar}The cyclic representation $\left(  S_{i}\right)
_{i=1}^{N}$ defined from the state $\omega$ in
\textup{(\ref{eqRemCo.0starstar})} has a co-invariant infinite-dimensional
subspace $\mathcal{L}$ spanned by the vectors $S_{I}^{\ast}\Omega$ where
$I=\left(  i_{1},\dots,i_{k}\right)  $ runs over all the multi-indices
$k\geq1$, and the relative orthocomplement $S\mathcal{L}\ominus\mathcal{L}$ is
spanned by the following $N^{2}$ orthogonal vectors:
\[
\left\{  S_{i}S_{j}^{\ast}\Omega:i,j=1,\dots,N\right\}  .
\]

\end{example}

\begin{proof}
The details are left to the reader.
\end{proof}

\begin{lemma}
\label{LemCo.1}Let $\left(  S_{i}\right)  _{i=1}^{N}$ be a representation of
$\mathcal{O}_{N}$ on a Hilbert space $\mathcal{H}$, and let $\mathcal{L}$ be a
closed subspace in $\mathcal{H}$ which is co-invariant. Then the subspace%
\begin{equation}
\mathcal{M}:=\bigvee\left\{  E_{i}\mathcal{L}:i=1,\dots,N\right\}
\label{eqCo.2}%
\end{equation}
is co-invariant and \emph{stable}. We use the notation $E_{i}=S_{i}S_{i}%
^{\ast}$ for the projection onto $S_{i}\mathcal{H}$.
\end{lemma}

\begin{proof}
Since $E_{i}^{2}=E_{i}$, it is clear that the space $\mathcal{M}$ is invariant
under each $E_{i}$. Since $S_{j}^{\ast}E_{i}=\delta_{j,i}S_{i}^{\ast}$, it is
clear that $\mathcal{M}$ is also invariant under each $S_{j}^{\ast}$. Hence
$\mathcal{M}$ is co-invariant. For every $x\in\mathcal{L}$, we have
$x=\sum_{i}E_{i}x$ by (\ref{eqCo.1}), and it follows then from (\ref{eqCo.2})
that $\mathcal{L}\subset\mathcal{M}$, which was used in the argument above.
\end{proof}

\section{\label{Res}Resolution subspaces of $\mathcal{H}$}

We begin this section with an explicit isomorphism between the family $%
\operatorname*{Co-inv}%
$ of all closed co-invariant subspaces, and the family $\operatorname*{Wan}$
of all closed wandering subspaces. The definitions refer to a specified
representation $\left\{  S_{i}\right\}  _{i=1}^{N}$ of $\mathcal{O}_{N}$ on a
Hilbert space $\mathcal{H}$, and we shall work with the corresponding
endomorphism $\alpha\left(  A\right)  =\sum_{i=1}^{N}S_{i}AS_{i}^{\ast}$,
$A\in\mathcal{B}\left(  \mathcal{H}\right)  $. Note that if $P$ is a
projection of $\mathcal{H}$ onto a closed subspace $\mathcal{K}\subset
\mathcal{H}$, then $\alpha\left(  P\right)  $ is the projection onto
$S\mathcal{K}$.

\begin{lemma}
\label{LemRes.0}For $\mathcal{L\in}%
\operatorname*{Co-inv}%
$ define $\mu\left(  \mathcal{L}\right)  =\left(  S\mathcal{L}\right)
\ominus\mathcal{L}$. Then $\mu\left(  \mathcal{L}\right)  \in
\operatorname*{Wan}$. For $\mathcal{W}\in\operatorname*{Wan}$, set
$\lambda\left(  \mathcal{W}\right)  =\mathcal{H}\ominus\mathcal{F}\left(
\mathcal{W}\right)  $; then $\lambda\left(  \mathcal{W}\right)  \in%
\operatorname*{Co-inv}%
$. Moreover%
\begin{equation}
\mathcal{W}=\mu\left(  \lambda\left(  \mathcal{W}\right)  \right)  .
\label{eqLemRes.0}
\end{equation}
\end{lemma}

\begin{proof}
If $\mathcal{K}\subset\mathcal{H}$ is a closed subspace, then we denote by
$P_{\mathcal{K}}$ the projection onto $\mathcal{K}$, i.e., $\mathcal{K}%
=P_{\mathcal{K}}\mathcal{H}$, and $P_{\mathcal{K}}=P_{\mathcal{K}}^{\ast
}=P_{\mathcal{K}}^{2}$. Thus, identifying closed subspaces with projections,
we arrive at the formulas
\[
\mu\left(  P_{\mathcal{L}}\right)  =\alpha\left(  P_{\mathcal{L}}\right)
-P_{\mathcal{L}}%
\]
and%
\[
\lambda\left(  P_{\mathcal{W}}\right)  =I_{\mathcal{H}}-P_{\mathcal{F}\left(
\mathcal{W}\right)  }=I_{\mathcal{H}}-\sum_{n=0}^{\infty}\alpha^{n}\left(
P_{\mathcal{W}}\right)  .
\]
If $\mathcal{W}$ is wandering, then the projections $\alpha^{n}\left(
P_{\mathcal{W}}\right)  $ in the sum are mutually orthogonal, and it follows
that $\lim_{n\rightarrow\infty}\alpha^{n}\left(  P_{\mathcal{W}}\right)  =0$.
Hence it follows that the two terms in the difference below are convergent,
and that%
\begin{equation}
\mu\left(  \lambda\left(  P_{\mathcal{W}}\right)  \right)  =\sum_{n=0}%
^{\infty}\alpha^{n}\left(  P_{\mathcal{W}}\right)  -\sum_{n=1}^{\infty}%
\alpha^{n}\left(  P_{\mathcal{W}}\right)  =P_{\mathcal{W}}.
\label{eqLemRes.0proof(3)}
\end{equation}

\end{proof}

\begin{lemma}
\label{LemRes.1}Let a representation $\left(  S_{i}\right)  $ of
$\mathcal{O}_{N}$ be given, and let $\mathcal{L}$ be a closed subspace which
is co-invariant. Let $\mathcal{V}$ denote the linear span of the spaces
$S^{k}\mathcal{L}$. Then $\mathcal{V}$ reduces to $\mathcal{L}$ in the representation.
\end{lemma}

\begin{proof}
We have $S^{\ast}\mathcal{L}\subset\mathcal{L}$. Using $\sum_{i}S_{i}%
S_{i}^{\ast}=I_{\mathcal{H}}$ we conclude that $\mathcal{L}\subset
S\mathcal{L}$, and by induction $S^{k}\mathcal{L}\subset S^{k+1}\mathcal{L}$
for $k\in\mathbb{N}$. Since, clearly, $S_{i}^{\ast}\left(  S^{k+1}%
\mathcal{L}\right)  \subset S^{k}\mathcal{L}$, we conclude that $\mathcal{V}$
reduces to $\mathcal{L}$ in the sense of Definition \ref{DefSub.1}%
(\ref{DefSub.1(7)}).
\end{proof}

The next two examples are the representation of $\mathcal{O}_{2}$ on the
Hilbert space $L^{2}\left(  \mathbb{T}\right)  $ which are defined from the
two best known examples of wavelets in $L^{2}\left(  \mathbb{R}\right)  $.

The first is the Haar wavelet. It has $\varphi=\chi_{I}^{{}}$, where $I$ is
the unit interval $I=\left[  0,1\right)  $, and%
\begin{equation}
\psi=\chi_{\left[  0,\frac{1}{2}\right)  }^{{}}-\chi_{\left[  \frac{1}%
{2},1\right)  }^{{}} \label{eqRes.1}%
\end{equation}
The corresponding functions $m_{j}$, $j=0,1$, which define the wavelet filters
are $m_{0}\left(  z\right)  =\frac{1}{\sqrt{2}}\left(  1+z\right)  $ and
$m_{1}\left(  z\right)  =\frac{1}{\sqrt{2}}\left(  1-z\right)  $; see
(\ref{eqInt.1}) and (\ref{eqInt.2}). The other wavelet is not localized in
$x$-space, but rather in the dual Fourier variable $\xi$ of the Fourier
transform $\hat{\psi}\left(  \xi\right)  =\int_{-\infty}^{\infty}e^{-i2\pi\xi
x}\psi\left(  x\right)  \,dx$. These wavelets are called frequency localized;
see \cite{BrJo99b}. It is known that there is a wavelet (named after Shannon)
for which the wavelet generator $\psi_{S}^{{}}$ is characterized by%
\begin{equation}
\hat{\psi}_{S}^{{}}\left(  \xi\right)  =\chi_{\left[  -\frac{1}{2},-\frac
{1}{4}\right)  \cup\left[  \frac{1}{4},\frac{1}{2}\right)  }^{{}}\left(
\xi\right)  ,\qquad\xi\in\mathbb{R}. \label{eqRes.2}%
\end{equation}
We now compare the $\mathcal{O}_{2}$-representation of $\psi_{S}^{{}}$ with
that of the Haar wavelet from (\ref{eqRes.1}).

More generally, we note that the examples of representations of $\mathcal{O}%
_{N}$ which primarily motivate our results are those which arise from
discretizing wavelet problems, in the sense of \cite{BrJo02b} and
\cite{Jor99b}. They may be realized on the Hilbert space $\mathcal{H}%
=L^{2}\left(  I\right)  $ where $I$ is a compact interval with Lebesgue
measure. Using the Fourier series on functions on $I$, it will be convenient
to view $\mathcal{H}$ alternately as $\ell^{2}\left(  \mathbb{Z}\right)  $, or
as $L^{2}\left(  \mathbb{T}\right)  $, where $\mathbb{T}$ is the one-torus,
equipped with the usual normalized Haar measure. When working with
$L^{2}\left(  \mathbb{T}\right)  $, it will be convenient to have the
orthonormal basis $e_{n}\left(  z\right)  =z^{n}$, $n\in\mathbb{Z}$, for use
in computations.

The following two examples of representations of $\mathcal{O}_{N}$ on
$\mathcal{H}$ illustrate the concepts in Section \ref{Int}. We will give them
just for $N=2$, but the reader can easily write out the general case.

\begin{example}
\label{ExaRes.2}Let $S_{0}f\left(  z\right)  =f\left(  z^{2}\right)  $, and
$S_{1}f\left(  z\right)  =zf\left(  z^{2}\right)  $, $f\in\mathcal{H}%
=L^{2}\left(  \mathbb{T}\right)  $.

In terms of the Fourier basis, this representation may be identified as
follows:%
\[
S_{0}e_{n}=e_{2n},\qquad S_{1}e_{n}=e_{2n+1},\qquad n\in\mathbb{Z},
\]
and with the adjoint operators%
\[
S_{0}^{\ast}e_{n}=\left\{
\begin{array}
[c]{ll}%
e_{n/2} & \text{if }n\text{ is even,}\\
0 & \text{if }n\text{ is odd,}%
\end{array}
\right.  \text{\qquad and\qquad}\left\{
\begin{array}
[c]{l}%
S_{1}^{\ast}e_{2n+1}=e_{n},\\
S_{1}^{\ast}e_{m}=0\text{ if }m\text{ is even.}%
\end{array}
\right.
\]
It is immediate that
the system of operators
$\left(  S_{i}\right)  _{i=0}^{1}$ defines a
representation of $\mathcal{O}_{2}$, and that $\mathcal{L}%
=\operatorname*{span}\left\{  e_{-1},e_{0}\right\}  $ is a \emph{core} for the
representation, in the sense of Definition \textup{\ref{DefSub.1}%
(\ref{DefSub.1(6)})}. Since%
\begin{equation}
\begin{aligned} S_{0}^{\ast}e_{-1} & =0, & \qquad S_{1}^{\ast}e_{-1} & =e_{-1},\\ S_{0}^{\ast}e_{0} & =e_{0}, & \qquad S_{1}^{\ast}e_{0} & =0, \end{aligned} \label{eqRes.3}%
\end{equation}
it is clear that $\mathcal{L}$ is \emph{stable} in the sense of Definition
\textup{\ref{DefSub.1}(\ref{DefSub.1(8)})}.
\end{example}

\begin{remark}
\label{RemRes.3}For Haar's representation $\left(  S_{i}\right)  _{i=0}^{1}$
in Example \textup{\ref{ExaRes.2}} we have two natural finite-dimensional
co-invariant subspaces,
\[
\mathcal{L}_{\min}=\bigvee\left\{  e_{-1}%
,e_{0}\right\}  \text{\qquad and\qquad }\mathcal{L}=\bigvee\left\{  e_{-2},e_{-1},e_{0}%
,e_{1}\right\}  ;
\]
and $\mathcal{L}_{\min}$ is minimal in the sense of
Definition \textup{\ref{DefSub.1}(\ref{DefSub.1(6)})}, while $\mathcal{L}$ is
not. This follows from the observations
\begin{align*}
\left(  S\mathcal{L}_{\min}\right)
\ominus\mathcal{L}_{\min} &= \bigvee\left\{  e_{-2},e_{1}\right\}  , \\
S\mathcal{L\ominus L} &= \bigvee\left\{  e_{-4},e_{-3},e_{2},e_{3}\right\}  , \\
\intertext{and}
S^{\ast}\mathcal{L} &= \mathcal{L}_{\min}.
\end{align*}
\end{remark}

\begin{example}
\label{ExaRes.4}Working with $I=\left[  -\frac{1}{2},\frac{1}{2}\right]  $, we
may view $L^{2}\left(  I\right)  $ as a space of $\mathbb{Z}$-periodic
functions on $\mathbb{R}$, and we will work with the operations $x\mapsto2x$,
$x\mapsto\frac{x}{2}$, $x\mapsto\frac{x+1}{2}$ modulo $\mathbb{Z}$. The last
two are branches of inverses of the doubling $x\mapsto2x\bmod1$, or
equivalently of $z\mapsto z^{2}$ on $\mathbb{T}$. In the second incarnation
the two branches of the inverse may be thought of as $\pm$ $\sqrt{z}$
restricted to $\mathbb{T}\subset\mathbb{C}$, i.e., the two branches of the
complex square root. Consider the subsets $A$ and $B$ of $I$: $A=\left\{  x\in
I:\left\vert x\right\vert \leq\frac{1}{4}\right\}  $, and $B=\left\{
x:-\frac{1}{2}\leq x<-\frac{1}{4}\right\}  \cup\left\{  x:\frac{1}{4}%
<x\leq\frac{1}{2}\right\}  $, forming a dyadic partition of $I$. Then set
\begin{equation}
\begin{aligned} S_{0}f\left( x\right) &:=\sqrt{2}\chi_{A}^{{}}\left( x\right) f\left( 2x\bmod{\mathbb{Z}}\right) \text{,\qquad and} \\ S_{1}f\left( x\right) &:=\sqrt{2}\chi_{B}^{{}}\left( x\right) f\left( 2x\bmod{\mathbb{Z}}\right). \end{aligned} \label{eqRes.4}%
\end{equation}
It will be convenient for us to omit the $\!{}\bmod{\mathbb{Z}}$ notation when
it is otherwise implicit in the formulas. Then it is immediate that the
isometries $\left(  S_{i}\right)  _{i=0}^{1}$ define a representation of
$\mathcal{O}_{2}$, and that $\mathcal{L}:=\mathbb{C}\openone$ is a \emph{core}
for the representation, acting on $\mathcal{H}=L^{2}\left(  I\right)  $, and
viewing $I$ as $\mathbb{R}\diagup\mathbb{Z}$. Also we use the notation
$\openone$ for the function which is constant${}\equiv1$ on $I$. The formulas%
\begin{align}
S_{0}^{\ast}\openone  &  =\frac{1}{\sqrt{2}}\openone, & S_{1}^{\ast}\openone
&  =\frac{1}{\sqrt{2}}\openone,\label{eqRes.5}\\
S_{0}\openone  &  =\sqrt{2}\chi_{A}^{{}}, & S_{1}\openone  &  =\sqrt{2}%
\chi_{B}^{{}} \label{eqRes.6}%
\end{align}
now make it clear that this one-dimensional core subspace $\mathcal{L}$ is not
\emph{stable} in the sense of Definition \textup{\ref{DefSub.1}%
(\ref{DefSub.1(8)})}. Specifically, we have:%
\begin{equation}
E_{0}\openone=\chi_{A}^{{}}\text{,\qquad and\qquad}E_{1}\openone=\chi_{B}^{{}%
}. \label{eqRes.7}%
\end{equation}
As in the discussion below Lemma \textup{\ref{LemRes.1}}, it follows that%
\begin{equation}
\begin{aligned} S_{0}^{\ast}\chi_{A}^{{}} & =\frac{1}{\sqrt{2}}\openone, & \qquad S_{0}^{\ast}\chi_{B}^{{}} & =0,\\ S_{1}^{\ast}\chi_{A}^{{}} & =0, & \qquad S_{1}^{\ast}\chi_{B}^{{}} & =\frac{1}{\sqrt{2}}\openone. \end{aligned} \label{eqRes.8}%
\end{equation}
Therefore the two-dimensional space $\mathcal{M}:=\operatorname*{span}\left\{
\chi_{A}^{{}},\chi_{B}^{{}}\right\}  $ is co-invariant, and it is
\emph{saturated} in the sense of Definition \textup{\ref{DefSub.1}%
(\ref{DefSub.1(5)})}. The space $\mathcal{M}$ is also stable.
\end{example}

\section{\label{Non}Nontrivial co-invariant subspaces}

Let $\left(  S_{i}\right)  _{i=1}^{N}$ be a representation of $\mathcal{O}%
_{N}$ acting on a Hilbert space $\mathcal{H}$. Suppose a closed subspace
$\mathcal{L}$ in $\mathcal{H}$ is co-invariant. We then get the resolution%
\begin{equation}
\mathcal{L}\subset S\mathcal{L}\subset S^{2}\mathcal{L}\subset\dots\subset
S^{n}\mathcal{L}\subset S^{n+1}\mathcal{L}\subset\dots\subset\mathcal{H}.
\label{eqNon.1}%
\end{equation}
If $\mathcal{L}$ is saturated, then $\bigvee_{n}S^{n}\mathcal{L}=\mathcal{H}$.

There are two trivial cases where both conditions are satisfied for a given
$\mathcal{L}$. First, if $\mathcal{L}=\mathcal{H}$, there isn't much to say.
Secondly, if $\mathcal{L}=S\mathcal{L}$, then the given representation on
$\mathcal{H}$ simply restricts to a representation on $\mathcal{L}$. Recall,
if $\mathcal{L}$ is co-invariant, and if $\mathcal{L}=S\mathcal{L}$, then
$\mathcal{L}$ is an invariant subspace for all the $2N$ operators
$\underbrace{S_{1},\dots,S_{N}}_{\text{isometries}}$, $\underbrace{S_{1}%
^{\ast},\dots,S_{N}^{\ast}}_{\text{co-isometries}}$. We say that $\mathcal{L}$
is a reducing subspace for the given representation. We then get two
orthogonal representations of $\mathcal{O}_{N}$ by restriction to each of the
two spaces, $\mathcal{L}$ and its orthocomplement $\mathcal{H}\ominus
\mathcal{L}=\left\{  x\in\mathcal{H}:\left\langle x\mid y\right\rangle
=0,\;y\in\mathcal{L}\right\}  $. We say that a co-invariant subspace
$\mathcal{L}$ for the given representation on $\mathcal{H}$ is
\emph{nontrivial} if $\mathcal{L}\neq\mathcal{H}$, and $\mathcal{L}\subsetneqq
S\mathcal{L}$.

\begin{theorem}
\label{ThmNon.1}There is a two-way correspondence between the nontrivial,
closed, co-invariant, and saturated subspaces $\mathcal{L}$ for a
representation $\left(  S_{i}\right)  _{i=1}^{N}$ of $\mathcal{O}_{N}$ on
$\mathcal{H}$, and closed subspaces $\mathcal{W}\neq0$ in $\mathcal{H}%
\ominus\mathcal{L}$ which provide an orthogonal decomposition%
\begin{equation}
\mathcal{H}=\mathcal{L}\oplus\mathcal{W}\oplus S\mathcal{W}\oplus
S^{2}\mathcal{W}\oplus\cdots\label{eqNon.2}%
\end{equation}
Note the subspace $\mathcal{W}$ in \textup{(\ref{eqNon.2})} is called
\emph{wandering}, and the conditions imply pairwise orthogonality of all the
closed subspaces:%
\begin{equation}
\mathcal{L},\;\mathcal{W},\;S_{i}\mathcal{W},\;S_{i_{1}}S_{i_{2}}%
\mathcal{W},\;S_{i_{1}}S_{i_{2}}S_{i_{3}}\mathcal{W},\;\dots,\;S_{i_{1}}\cdots
S_{i_{k}}\mathcal{W},\;\dots, \label{eqNon.3}%
\end{equation}
where all multi-indices $\left(  i_{1},i_{2},\dots,i_{k}\right)  $ are
considered for $k=1,2,\dots$.
\end{theorem}

\begin{proof}
Suppose first that $0\neq\mathcal{W}$ satisfies the conditions from
(\ref{eqNon.2}). Then set $\mathcal{L}:=\mathcal{H}\ominus
\sideset{}{^{\smash{\oplus}}}{\textstyle\sum}\limits_{k=0}^{\infty}%
S^{k}\mathcal{W}$, with the understanding that $S^{0}\mathcal{W}:=\mathcal{W}$
and
\begin{equation}
S^{k}\mathcal{W}=\bigvee_{\left(  i_{1},i_{2},\dots,i_{k}\right)  }S_{i_{1}%
}\cdots S_{i_{k}}\mathcal{W}. \label{eqNon.4}%
\end{equation}
Since $\sideset{}{^{\smash{\oplus}}}{\textstyle\sum}\limits_{k=0}^{\infty
}S^{k}\mathcal{W}$ is invariant under all the isometries $S_{i}$,
$i=1,\dots,N$, it is clear that $\mathcal{L}$ is invariant under the adjoints
$S_{i}^{\ast}$, and so $\mathcal{L}$ is co-invariant. Since $\mathcal{W}\neq
0$, clearly $\mathcal{L}\neq\mathcal{H}$. It remains to verify that the strict
inclusion $\mathcal{L}\subsetneqq S\mathcal{L}$ holds. Suppose, indirectly,
that $\mathcal{L}=S\mathcal{L}$. Then $\mathcal{L}$ is reducing, and so is
$\sideset{}{^{\smash{\oplus}}}{\textstyle\sum}\limits_{k=0}^{\infty}%
S^{k}\mathcal{W}$. To see that this is impossible, let $\mathcal{I}\left(
k\right)  $ denote the set of all multi-indices $\left(  i_{1},\dots
,i_{k}\right)  $, where $i_{1},i_{2},\dots\in\left\{  1,2,\dots,N\right\}  $.
If $\sideset{}{^{\smash{\oplus}}}{\textstyle\sum}\limits_{k=0}^{\infty}%
S^{k}\mathcal{W}$ were $S_{i}^{\ast}$-invariant for all $i$, then
\begin{equation}
\lim_{k\rightarrow\infty}\sum_{I\in\mathcal{I}\left(  k\right)  }\left\Vert
S_{I}^{\ast}x\right\Vert ^{2}=0\text{\qquad for all }x\in\mathcal{W},
\label{eqNon.5}%
\end{equation}
where $S_{I}^{\ast}=S_{i_{k}}^{\ast}\cdots S_{i_{1}}^{\ast}$ for $I=\left(
i_{1},\dots,i_{k}\right)  $. But this is impossible by the Cuntz relations.
Recall that
\begin{equation}
\sum_{I\in\mathcal{I}\left(  k\right)  }\left\Vert S_{I}^{\ast}x\right\Vert
^{2}=\left\Vert x\right\Vert ^{2}\text{\qquad holds for all }x.
\label{eqNon.6}%
\end{equation}
This contradicts (\ref{eqNon.5}) and the condition $\mathcal{W}\neq0$.

Suppose now that $\mathcal{L}$ is a given nontrivial, closed, co-invariant and
saturated subspace for the given representation in $\mathcal{H}$. As noted,
then $\mathcal{L}\subsetneqq S\mathcal{L}$; so the relative orthocomplement
$\mathcal{W}:=\left(  S\mathcal{L}\right)  \ominus\mathcal{L}$ is nonzero. We
first prove that all the multi-indexed subspaces listed in (\ref{eqNon.3}) are
mutually orthogonal. The argument is by induction, starting with
$\mathcal{W}\perp S\mathcal{W}$. We must check the inner products%
\begin{equation}
\ip{\sum_{i}S_{i}x_{i}}{\sum_{j_{1},j_{2}}S_{j_{1}}S_{j_{2}}y_{j_{1},j_{2}}}
\label{eqNon.7}%
\end{equation}
where $x_{i},y_{j_{1},j_{2}}\in\mathcal{L}$, under the conditions%
\begin{equation}
\sum_{i}S_{i}x_{i}\perp\mathcal{W}\text{\qquad and\qquad}\sum_{j_{2}}S_{j_{2}%
}y_{j_{1},j_{2}}\perp\mathcal{W}\text{\quad for all }j_{1}=1,\dots,N.
\label{eqNon.8}%
\end{equation}
But the term in (\ref{eqNon.7}) simplifies to%
\[
\sum_{i,j}\left\langle x_{i}\mid S_{j}y_{i,j}\right\rangle =0,
\]
where the vanishing results from the conditions (\ref{eqNon.8}). We leave the
remaining recursive argument to the reader.

It remains to check that the sequence of closed subspaces%
\[
\mathcal{L},\;\mathcal{W},\;S\mathcal{W},\;S^{2}\mathcal{W},\;\dots
\]
is total in $\mathcal{H}$, i.e., that (\ref{eqNon.2}) holds.

Since $\mathcal{L}$ is given to be saturated, the conclusion will follow from
\begin{equation}
S^{k+1}\mathcal{L}\subset\mathcal{L}\oplus\mathcal{W}\oplus S\mathcal{W}%
\oplus\dots\oplus S^{k}\mathcal{W}. \label{eqNon.9}%
\end{equation}
If $k=0$, this holds from the ansatz which defines $\mathcal{W}$ in terms of
the given closed subspace $\mathcal{L}$, i.e., $\mathcal{W}:=\left(
S\mathcal{L}\right)  \ominus\mathcal{L}$. Suppose (\ref{eqNon.9}) has been
verified up to $k-1$. Then%
\begin{align*}
S^{k+1}\mathcal{L}  &  \subset S\left(  \mathcal{L}\oplus\mathcal{W}%
\oplus\dots\oplus S^{k-1}\mathcal{W}\right) \\
&  \subset
\smash{\underset{{}\rlap{$\scriptstyle\searrow$}}{S\mathcal{L}}}\oplus
S\mathcal{W}\oplus\dots\oplus S^{k}\mathcal{W}\\
&  \subset\overbrace{\mathcal{L}\oplus\mathcal{W}}\oplus S\mathcal{W}%
\oplus\dots\oplus S^{k}\mathcal{W},
\end{align*}
where we used the induction hypothesis and the ansatz.
\end{proof}

\section{\label{Exi}Existence}

While the co-invariant subspaces for the wavelet representations are
relatively well understood for the special representations of $\mathcal{O}%
_{N}$ which derive from wavelet analysis, the situation is somewhat mysterious
in the case of the most general $\mathcal{O}_{N}$-representations. Even the
existence of nontrivial co-invariant subspaces which form cores for a given
$\mathcal{O}_{N}$-representation appears to be open in general. However:

\begin{corollary}
\label{CorExi.1}A representation $\left(  S_{i}\right)  _{i=1}^{N}$ of
$\mathcal{O}_{N}$ on a Hilbert space $\mathcal{H}$ has a nontrivial \emph{core
subspace} $\mathcal{L}$, i.e., a closed co-invariant subspace $\mathcal{L}%
\subset\mathcal{H}$ which is nontrivial and saturated,
if and only if it has a nonzero wandering vector. Further, $\mathcal{L}$
may be chosen to be minimal with respect to being co-invariant and saturated.
\end{corollary}

\begin{proof}
We showed in Theorem \ref{ThmNon.1} that the general case may be reduced to
the consideration of nontrivial co-invariant subspaces. In Theorem
\ref{ThmNon.1}, we also showed that the nontrivial co-invariant subspaces
$\mathcal{L}$ may be understood from the corresponding wandering subspaces
$\mathcal{W}$. If $\mathcal{L}$ is a nontrivial co-invariant subspace, then
$\mathcal{W}:=\left(  S\mathcal{L}\right)  \ominus\mathcal{L}$ is wandering
for the representation, and conversely every wandering subspace $\mathcal{W}$
comes from a co-invariant subspace $\mathcal{L}$. The idea is to prove the
existence claim in the corollary by establishing the existence of a maximal
wandering subspace $\mathcal{W}_{\max}$. We will do this by
Zorn's lemma. We say
that a wandering subspace $\mathcal{W}_{\max}$ is maximal if for every
wandering subspace $\mathcal{W}$ such that $\mathcal{W}\supseteq
\mathcal{W}_{\max}$, we may conclude that $\mathcal{W}=\mathcal{W}_{\max}$.
The existence of $\mathcal{W}_{\max}$ follows by an application of Zorn's
lemma to the family of all wandering subspaces $\mathcal{W}\neq0$ ordered by
inclusion. To be able to do this, we need to know that every representation of
$\mathcal{O}_{N}$ has at least one wandering subspace $\mathcal{W}\neq0$.
Some details are given below for the benefit of the reader.

To start the transfinite induction we must assume that there are solutions%
\[
x\in\mathcal{H},\qquad\left\Vert x\right\Vert =1,
\]
to the system%
\begin{equation}
\left\langle x\mid S_{I}x\right\rangle =0\text{,\qquad for all }I=\left(
i_{1},\dots,i_{k}\right)  ,\;k\geq1. \label{eqExi.starstarstar}%
\end{equation}
This is the stated condition in the 
formulation of Corollary \ref{CorExi.1}, and it is needed for the start of
the transfinite induction. It is conceivable that such a nonzero vector
$ x $
automatically exists, but we have not been abe to show this. We thank
Ken Davidson for calling this point to our attention. 
\end{proof}

\begin{remark}
\label{RemExiFeb18}It follows from Section \ref{Co} above that
$\left\Vert 1-s_{I}\right\Vert\geq\sqrt{2}$
for all $I$ where $s_{I}\in\mathcal{O}_{N}$ satisfies
$\rho\left(  s_{I}\right)  =S_{I}=S_{i_{1}}\cdots S_{i_{k}}$. Hence it follows
that there is a state $\omega$ on $\mathcal{O}_{N}$ such that
$\omega\left( S_{I}\right) =0$
for all $I\in\mathcal{I}\left( N\right) \setminus\left\{ \varnothing\right\} $,
but it is not known if
such a normal state $\omega$ may be found. If so
we would have
\[
\operatorname*{trace}\left(  \left\vert x\right\rangle \left\langle
x\right\vert S_{I}\right)  =\left\langle x\mid S_{I}x\right\rangle 
=0,
\]
where $\left\vert x\right\rangle \left\langle x\right\vert $ is the rank-$1$
projection of $\mathcal{H}$ onto $\mathbb{C}x$, and where we used Dirac's
bra-ket notation.
\end{remark}

\begin{remark}
\label{RemExiJan29.end}The existence problem
for wandering subspaces $\mathcal{W}$ may be
reformulated as a fixed-point
problem. This is made clear by the
equivalent identities \textup{(\ref{eqLemRes.0})} and
\textup{(\ref{eqLemRes.0proof(3)})} of Lemma \textup{\ref{LemRes.0}}.
Introducing the transformation
$\mathcal{W}\mapsto\mu\left( \lambda\left( \mathcal{W}\right) \right) $ of
the right-hand side in \textup{(\ref{eqLemRes.0})}, or
$P_{\mathcal{W}}\mapsto\mu\left( \lambda\left( P_{\mathcal{W}}\right) \right) $
in \textup{(\ref{eqLemRes.0proof(3)})}, it is clear
that the existence of a wandering subspace
$\mathcal{W}$ is equivalent to the existence of
a fixed point $\neq0$ in the identity
\textup{(\ref{eqLemRes.0})}.
\end{remark}

\section{\label{Pur}Pure co-invariant subspaces}

Let $\left(  S_{i}\right)  _{i=1}^{N}$ be a representation of $\mathcal{O}%
_{N}$ on a Hilbert space $\mathcal{H}$, and let $\mathcal{L}$ be a
co-invariant subspace. For $x\in\mathcal{H}$, set%
\begin{align*}
S^{n}x  &  :=\bigvee\left[  S_{i_{1}}\cdots S_{i_{n}}x;I=\left(  i_{1}%
,\dots,i_{k}\right)  \in\mathcal{I}_{n}\left(  N\right)  \right]  ,\\
S^{0}x  &  :=x,\\%
\intertext{and}%
S^{\infty}x  &  :=\bigvee_{n\geq0}S^{n}x.
\end{align*}
The special case $N=1$ covers a single isometry.

\begin{definition}
\label{DefPur.1}We say that a co-invariant subspace $\mathcal{L}$ is
\emph{pure} if the following implication holds:%
\begin{equation}
x\in\mathcal{L},\;S^{\infty}x\subset\mathcal{L}\Longrightarrow x=0.
\label{eqPur.1}%
\end{equation}

\end{definition}

To motivate the next result, we recall first the case $N=1$, i.e., a single
isometry $S$ in a Hilbert space $\mathcal{H}$. The Wold decomposition states
that $\mathcal{H}$ decomposes as $\mathcal{H}=\mathcal{H}_{0}\oplus
\mathcal{H}_{\infty}$, $\mathcal{H}_{\infty}=\bigwedge_{n\geq1}S^{n}%
\mathcal{H}$, $\mathcal{H}_{0}=\mathcal{H}\ominus\mathcal{H}_{\infty}$, both
subspaces invariant, with $S|_{\mathcal{H}_{\infty}}$ unitary, and
$S_{0}:=S|_{\mathcal{H}_{0}}$ a shift. To say that $S_{0}$ is a \emph{shift}
means $S_{0}^{\ast\,n}\underset{n\rightarrow\infty}{\longrightarrow}0$. A
shift is determined up to unitary equivalence by its multiplicity, i.e.,
$\dim\left\{  x\in\mathcal{H}_{0};S_{0}^{\ast}x=0\right\}  $. The space
$\mathcal{W}_{0}:=\ker\left(  S_{0}^{\ast}\right)  =\mathcal{H}_{0}\ominus
S_{0}\mathcal{H}_{0}$ is wandering, and $\left(  \mathcal{H}_{0},S_{0}\right)
$ is unitarily equivalent to the obvious shift on vectors, i.e., to
$\sideset{}{^{\oplus}}{\textstyle\sum}\limits_{n\geq0}\mathcal{W}%
_{0}:=\mathcal{W}_{0}\oplus\mathcal{W}_{0}\oplus\cdots$ and the operator given
by $\left(  x_{0},x_{1},\dots\right)  \mapsto\left(  0,x_{0},x_{1}%
,\dots\right)  $. So for $N=1$, our understanding of the invariant and the
co-invariant subspaces amounts to the corresponding issues for the shift.
Since a subspace $\mathcal{L}$ is co-invariant if and only if its
orthocomplement is invariant, the problem is solved by Beurling's theorem.
Identify $\smash[b]{\sideset{}{^{\oplus}}{\textstyle\sum}\limits_{n\geq
0}}\mathcal{W}_{0}$ by the Hardy space $H_{+}\left(  \mathcal{W}_{0}\right)  $
of analytic $\mathcal{W}_{0}$-valued functions, i.e.,
\begin{equation}
f\colon\mathbb{T}\longrightarrow\mathcal{W}_{0},\;f\left(  z\right)
=\sum_{n=0}^{\infty}\xi_{n}z^{n},\qquad\xi_{n}\in\mathcal{W}_{0},\;\left\Vert
f\right\Vert ^{2}=\sum_{n=0}^{\infty}\left\Vert \xi_{n}\right\Vert ^{2}.
\label{eqPur.2}%
\end{equation}
An inner function $u$ is a function on $\mathbb{T}$, taking values in the
unitary operators on $\mathcal{W}_{0}$, i.e., $u\left(  z\right)  ^{\ast
}u\left(  z\right)  =I_{\mathcal{W}_{0}}$, $\mathrm{a.e.}\;z\in\mathbb{T}$,
with $u$ having an analytic operator-valued continuation to $\left\{
z\in\mathbb{C};\left\vert z\right\vert <1\right\}  $. Since, by (\ref{eqNon.2}%
), the shift is represented by
\[
M_{z}f\left(  z\right)  =zf\left(  z\right)
, \qquad f\in H_{+}\left(  \mathcal{W}_{0}\right)  ,
\]
it follows that the space
\[
H_{+}\left(  u\right)  :=\left\{  uf;f\in H_{+}\left(  \mathcal{W}%
_{0}\right)  \right\}  
\]
is invariant when $u$ is an inner function.
Beurling's theorem \cite{Hel64} states that every invariant subspace for the
shift has this form.

\begin{proposition}
\label{ProPur.2}If an inner function $u$ is given, then the co-invariant
subspace%
\begin{equation}
\mathcal{L}\left(  u\right)  :=H_{+}\left(  u\right)  ^{\perp}=H_{+}\left(
\mathcal{W}_{0}\right)  \ominus H_{+}\left(  u\right)  \label{eqPur.3}%
\end{equation}
is pure.
\end{proposition}

\begin{proof}
If $u$ is constant, $u\equiv I_{\mathcal{W}_{0}}$, then $\mathcal{L}\left(
u\right)  =0$ satisfies the condition, so we assume $u$ to be non-constant.
Suppose $f\in\mathcal{L}\left(  u\right)  $ and $S^{\infty}f\subset
\mathcal{L}\left(  u\right)  $. Then%
\begin{equation}
f\left(  z\right)  ,\;zf\left(  z\right)  ,\;z^{2}f\left(  z\right)
,\;\dots\text{\qquad are all in }\mathcal{L}\left(  u\right)  .
\label{eqPur.4}%
\end{equation}
Let $f=\sum_{n=0}^{\infty}\xi_{n}z^{n}$, $\xi_{n}\in\mathcal{W}_{0}$, and
$u\left(  z\right)  =\sum_{n=0}^{\infty}A_{n}z^{n}$ where $A_{n}%
\colon\mathcal{W}_{0}\rightarrow\mathcal{W}_{0}$ is a system of operators in
$\mathcal{W}_{0}$ such that%
\begin{align}
\sum_{n=0}^{\infty}A_{n}A_{n}^{\ast}  &  =I_{\mathcal{W}_{0}}\text{\qquad
and}\label{eqPur.5}\\
\sum_{n=0}^{\infty}A_{n}A_{n+k}^{\ast}  &  =0\text{\qquad for }k=1,2,\dots.
\label{eqPur.6}%
\end{align}
Using (\ref{eqPur.3}) and (\ref{eqPur.4}), we get the system:%
\begin{align*}
A_{n}^{\ast}\xi_{0}+A_{n+1}^{\ast}\xi_{1}+\dots &  =0,\\
A_{n}^{\ast}\xi_{1}+A_{n+1}^{\ast}\xi_{2}+\dots &  =0,\\
A_{n}^{\ast}\xi_{2}+A_{n+1}^{\ast}\xi_{3}+\dots &  =0,\\
&  \vdots\;.
\end{align*}
Multiply each of these equations by $A_{n}$ and sum over $n=0,1,2,\dots$. We
then get%
\[
\sum_{n=0}^{\infty}A_{n}A_{n}^{\ast}\xi_{k}+\sum_{n=0}^{\infty}A_{n}%
A_{n+1}^{\ast}\xi_{k+1}+\dots=0
\]
for $k=0,1,\dots$, An application of (\ref{eqPur.5})--(\ref{eqPur.6}) now
yields $\xi_{k}=0$ for all $k=0,1,\dots$. This proves the conclusion.
\end{proof}

If $N>1$, we are considering representations of $\mathcal{O}_{N}$, and the
examples in Section \ref{Sub} show that not all co-invariant subspaces are
pure. But we do have the following:

\begin{theorem}
\label{ThmPur.3}Let $\left(  S_{i}\right)  _{i=1}^{N}$, $N>1$, be a
representation of $\mathcal{O}_{N}$ on a Hilbert space $\mathcal{H}$, and let
$\mathcal{L}$ be a co-invariant subspace. Then $\mathcal{L}$ contains a pure
co-invariant subspace.
\end{theorem}

\begin{proof}
If $\mathcal{L}=0$, we are done. If not, consider vectors (if any)
$x\in\mathcal{L}$, $\left\Vert x\right\Vert =1$, such that $S^{\infty}%
x\subset\mathcal{L}$. If there are no such vectors, $\mathcal{L}$ is pure.
Otherwise pick a family $\mathcal{F}=\left\{  x\right\}  $ such that
$S^{\infty}x\subset\mathcal{L}$ and $S^{\infty}x\perp S^{\infty}y$ for all
$x,y\in\mathcal{F}$, $x\neq y$. By Zorn's lemma, we may pick $\mathcal{F}$ to
be maximal with respect to these properties. Now set $\mathcal{L}%
_{p}:=\mathcal{L}\ominus\left[  S^{\infty}x;x\in\mathcal{F}\right]  $. This
space is clearly co-invariant. If some vector $x_{p}\in\mathcal{L}_{p}$
satisfies $S^{\infty}x_{p}\subset\mathcal{L}_{p}$, then $\mathcal{F}%
\cup\left\{  x_{p}\right\}  $ satisfies the orthogonality property. Since
$\mathcal{F}$ was chosen maximal, we conclude that $x_{p}=0$. This proves that
$\mathcal{L}_{p}$ is pure.
\end{proof}

\begin{examples}
\label{ExaPur.4} \textup{(a)}~Consider the representation $\left(
S_{i}\right)  _{i=0}^{1}$ of $\mathcal{O}_{2}$ acting on $\ell^{2}\left(
\mathbb{Z}\right)  \simeq L^{2}\left(  \mathbb{T}\right)  $ outlined in
Example \textup{\ref{ExaRes.2}}. From \textup{(\ref{eqRes.3})} we note that
the closed subspace spanned by $\left\{  e_{0},e_{1},e_{2},\dots\right\}  $ is
invariant for the representation, i.e., invariant under all four of the
operators $S_{i}$, $S_{j}^{\ast}$, $i,j=0,1$; and moreover this restricted
representation, $\rho_{+}^{{}}$ say, is irreducible on this space. The space
is $H_{+}\simeq\ell^{2}\left(  \left\{  0,1,2,\dots\right\}  \right)  \simeq
{}$the Hardy space. Recall the Hardy space $H_{+}\subset L^{2}\left(
\mathbb{T}\right)  $ has the representation%
\begin{equation}
H_{+}=\left\{  f\in L^{2}\left(  \mathbb{T}\right)  :f\left(  z\right)
=\sum_{n=0}^{\infty}\xi_{n}z^{n},\;z\in\mathbb{T},\;\sum_{n=0}^{\infty
}\left\vert \xi_{n}\right\vert ^{2}\;(=\left\Vert f\right\Vert ^{2}%
)\;<\infty\right\}  . \label{eqPur.7}%
\end{equation}
Since, in general, \emph{every finite-dimensional} co-invariant subspace is
pure, we note that each one of the following subspaces $\mathcal{L}_{1}$,
$\mathcal{L}_{2}$, and $\mathcal{L}_{3}$ in $H_{+}$ is pure: $\mathcal{L}%
_{1}=\left[  e_{0}\right]  $, $\mathcal{L}_{2}=\left[  e_{0},e_{1}\right]  $,
$\mathcal{L}_{3}=\left[  e_{0},e_{1},e_{2}\right]  $. The corresponding three
wandering subspaces, i.e.,
\begin{equation}
\mathcal{W}_{i}:=\left(  S\mathcal{L}_{i}\right)  \ominus\mathcal{L}%
_{i},\qquad i=1,2,3, \label{eqPur.8}%
\end{equation}
are:%
\begin{equation}
\mathcal{W}_{1}=\left[  e_{1}\right]  ,\qquad\mathcal{W}_{2}=\left[
e_{2},e_{3}\right]  ,\qquad\mathcal{W}_{3}=\left[  e_{3},e_{4},e_{5}\right]  .
\label{eqPur.8bis}%
\end{equation}
An example of a proper \textup{(}i.e., nontrivial\/\textup{)} co-invariant
subspace ($\subset H_{+}$) which is not pure is%
\begin{equation}
\mathcal{L}:=\left[  e_{0},e_{1},e_{3},e_{6},e_{7},e_{12},e_{13},e_{14}%
,e_{15},e_{24},e_{25},\dots\right]  . \label{eqPur.9}%
\end{equation}
The stated properties for this last subspace $\mathcal{L}$ follow from its
representation as%
\begin{equation}
\mathcal{L}=\mathcal{L}_{2}\oplus S^{\infty}e_{3}. \label{eqPur.10}%
\end{equation}

\textup{(b)}~Based on \textup{(a)}, one might think that a pure co-invariant
subspace cannot be infinite-dimensional. This is not so, as we now illustrate
with the representation $\rho$ from Remark \textup{\ref{RemCo.0}}. This is the
representation$\ \rho$ of $\mathcal{O}_{N}$ with cyclic vector $\Omega$ which
acts on the Hilbert space $\mathcal{H}$, i.e., $\mathcal{H}$ is the Hilbert
space spanned by the following vectors: $S_{I}S_{J}^{\ast}\Omega$, where $I$
and $J$ run over the set of all multi-indices $\mathcal{I}\left(  N\right)  $,
i.e., $\left(  i_{1},i_{2},\dots,i_{n}\right)  $, $n=1,2,\dots$, $i_{\nu}%
\in\left\{  1,2,\dots,N\right\}  $, with the convention $n=0$ corresponding to
$I=\varnothing$, and $S_{\varnothing}=I$. The two indices $I$ and $J$ might
have different length. It is immediate from \textup{(\ref{eqRemCo.0starstar})}
that the following subspace,
\begin{equation}
\mathcal{L}:=\bigvee\left[  S_{I}^{\ast}\Omega:I\in\mathcal{I}\left(
N\right)  \setminus\left\{  \varnothing\right\}  \right]  , \label{eqPur.11}%
\end{equation}
is a nontrivial co-invariant subspace for the representation $\rho$, and that
the corresponding wandering subspace $\mathcal{W}$ is $N^{2}$-dimensional; in
fact, $\mathcal{W}=\bigvee_{i,j}\left[  S_{i}S_{j}^{\ast}\Omega\right]  $; see
Example \ref{ExaCo.0starstarstar}.
\end{examples}

\begin{observation}
\label{ObsPur.5}The subspace $\mathcal{L}$ in \textup{(\ref{eqPur.11})} is a
\emph{pure} co-invariant subspace for the representation $\rho$ of Remark
\textup{\ref{RemCo.0}}.
\end{observation}

\begin{proof}
Let $P=P_{\mathcal{L}}$ denote the (orthogonal) projection onto $\mathcal{L}$.
We show that%
\begin{equation}
\left\Vert PS_{i}P\right\Vert \leq N^{-\frac{1}{2}}\text{\qquad for }%
i=1,\dots,N. \label{eqPur.12}%
\end{equation}
Since $N>1$, the result follows; in fact, if $x$ is any vector $x\in
\mathcal{L}$ for which $S_{i}x\in\mathcal{L}$ for some $i$, then $x=0$.
Indeed, $\left\Vert x\right\Vert =\left\Vert S_{i}x\right\Vert =\left\Vert
PS_{i}Px\right\Vert \leq N^{-\frac{1}{2}}\left\Vert x\right\Vert $.

To prove (\ref{eqPur.12}), we note that the following normalized vectors,%
\begin{equation}
\left\{  N^{\frac{k}{2}}S_{I}^{\ast}\Omega:\left\vert I\right\vert =k^{\prime
}\right\}  , \label{eqPur.13}%
\end{equation}
are mutually orthogonal when $I$ varies over $\mathcal{I}_{k^{\prime}}\left(
N\right)  $, and also when the respective length of indices $I$, $I^{\prime}$
are different, i.e., $\left\vert I\right\vert =k\neq k^{\prime}=\left\vert
I^{\prime}\right\vert $. Hence, by (\ref{eqPur.11}), we have an orthonormal
basis for $\mathcal{L}$. Now let $P_{k}$ denote the projection onto the closed
subspace spanned by the vectors in (\ref{eqPur.13}) for all values $k^{\prime
}$ such that $k^{\prime}\leq k$. Then we get
\begin{equation}
\lim_{k\rightarrow\infty}P_{k}=P_{\mathcal{L}}\;(=P). \label{eqPur.14}%
\end{equation}
If $x\in\mathcal{L}$, then
\begin{align*}
\left\Vert P_{k}S_{i}x\right\Vert ^{2}  &  =N^{k}\sum_{\left\vert I\right\vert
=k}\left\vert \left\langle S_{I}^{\ast}\Omega\mid S_{i}x\right\rangle
\right\vert ^{2}\\
&  =N^{-1}N^{k+1}\sum_{\left\vert I\right\vert =k}\left\vert \left\langle
\smash{\underbrace{S_{i_{\mathstrut}}^{\ast }S_{I}}_{\scriptscriptstyle\!\!\!\left( k+1\right) \text{-index}\!\!\!}}\Omega
\mid x\right\rangle \right\vert ^{2}%
\vphantom{\underbrace{S_{i_{\mathstrut}}^{\ast }S_{I}}_{\scriptscriptstyle\!\!\!\left( k+1\right) \text{-index}\!\!\!}}\\
&  \leq N^{-1}\left\Vert P_{k+1}x\right\Vert ^{2}.
\end{align*}
Letting $k\rightarrow\infty$, and using (\ref{eqPur.14}), we now arrive at the
conclusions: $x=P_{\mathcal{L}}x$ ($=\lim_{k\rightarrow\infty}P_{k}x$),
$P_{\mathcal{L}}S_{i}x=\lim_{k\rightarrow\infty}P_{k}S_{i}x$, and $\left\Vert
P_{\mathcal{L}}S_{i}x\right\Vert \leq N^{-\frac{1}{2}}\left\Vert x\right\Vert
$, the last estimate being equivalent to the desired one (\ref{eqPur.12}).
\end{proof}

\section{\label{Tig}Tight frames of wavelets}

We now turn to the representations (\ref{eqInt.2})--(\ref{eqInt.4}) which
define tight frames of multiresolution wavelets, and we give a representation
of pure co-invariant subspaces in the Hilbert space $L^{2}\left(
\mathbb{T}\right)  $. It is shown in \cite{BrJo02b} that to get solutions
$\varphi$, $\psi_{i}^{{}}$ as in (\ref{eqInt.1}) which are in $L^{2}\left(
\mathbb{R}\right)  $, the following condition must be satisfied by the matrix
function $\mathbb{T}\ni z\mapsto A\left(  z\right)  \in\mathrm{U}_{N}\left(
\mathbb{C}\right)  $ from formula (\ref{eqInt.4}). Specifically, if $\rho
=\rho_{N}^{{}}=e^{i2\pi/N}$, $i=\sqrt{-1}$, and if%
\begin{equation}
A_{j,k}\left(  1\right)  =\frac{1}{\sqrt{N}}\rho^{j\cdot k}=\frac{1}{\sqrt{N}%
}e^{i\frac{j\cdot k2\pi}{N}}, \label{eqTig.1}%
\end{equation}
then the solutions $\hat{\varphi}$, $\hat{\psi}_{i}^{{}}$ in (\ref{eqInt.4})
are in $L^{2}\left(  \mathbb{R}\right)  $, and their inverse Fourier
transforms $\varphi$, $\psi_{i}^{{}}$ are in $L^{2}\left(  \mathbb{R}\right)
$ as well. Introducing the triple-indexed functions%
\begin{equation}
\psi_{i,j,k}^{{}}\left(  x\right)  =N^{\frac{j}{2}}\psi_{i}^{{}}\left(
N^{j}x-k\right)  ,\qquad i=1,\dots,N-1,\;j\in\mathbb{Z},\;k\in\mathbb{Z},
\label{eqTig.2}%
\end{equation}
we say that $\left\{  \psi_{i,j,k}\right\}  $ is a \emph{tight frame} if and
only if, for all $F\in L^{2}\left(  \mathbb{R}\right)  $, the following
(Bessel) identity holds:%
\begin{equation}
\int_{\mathbb{R}}\left\vert F\left(  x\right)  \right\vert ^{2}\,dx=\sum
_{i,j,k}\left\vert \left\langle \psi_{i,j,k}^{{}}\mid F\right\rangle
_{L^{2}\left(  \mathbb{R}\right)  }\right\vert ^{2} \label{eqTig.3}%
\end{equation}
The next result is also proved in \cite{BrJo02b}, but is included here for the
convenience of the reader:

\begin{theorem}
\label{ThmTig.1}If $\mathbb{T}\ni z\mapsto A\left(  z\right)  \in
\mathrm{U}_{N}\left(  \mathbb{C}\right)  $ is a unitary matrix function
satisfying \textup{(\ref{eqTig.1})} then the wavelet functions $\psi
_{i,j,k}^{{}}$ defined in \textup{(\ref{eqInt.1})} and \textup{(\ref{eqTig.2}%
)} form a tight frame in $L^{2}\left(  \mathbb{R}\right)  $.
\end{theorem}

\begin{remark}
\label{RemTig.2}\textup{(The stretched Haar wavelet.)} The following example
for $N=2$ shows that this system \textup{(\ref{eqTig.2})} might not in fact be
an orthonormal basis: Take $A\left(  z\right)  =\dfrac{1}{\sqrt{2}}%
\begin{pmatrix}
1 & z\\
1 & -z
\end{pmatrix}
$. Then $A\left(  1\right)  =\dfrac{1}{\sqrt{2}}%
\begin{pmatrix}
1 & 1\\
1 & -1
\end{pmatrix}
$,
\begin{align}
\varphi\left(  x\right)   &  =\frac{1}{3}\chi_{\left[  0,3\right)  }^{{}%
}\left(  x\right)  \text{,\qquad and}\label{eqTig.4}\\
\psi\left(  x\right)   &  =\frac{1}{3}\left(  \chi_{\left[  0,\frac{3}%
{2}\right)  }^{{}}\left(  x\right)  -\chi_{\left[  \frac{3}{2},3\right)  }%
^{{}}\left(  x\right)  \right)  . \label{eqTig.5}%
\end{align}
It follows from a direct verification, or from the theorem, that this function
$\psi$ in \textup{(\ref{eqTig.5})} makes%
\begin{equation}
\psi_{j,k}^{{}}\left(  x\right)  :=2^{\frac{j}{2}}\psi\left(  2^{j}x-k\right)
,\qquad j,k\in\mathbb{Z}, \label{eqTig.6}%
\end{equation}
into a tight frame. But since $\left\Vert \psi\right\Vert =\frac{1}{\sqrt{3}}%
$, and since the different functions in \textup{(\ref{eqTig.6})} are
\emph{not} orthogonal, we see that this is a wavelet tight frame which is not
an orthonormal basis in $L^{2}\left(  \mathbb{R}\right)  $.
\end{remark}

Not all representations as in (\ref{eqInt.2}) and (\ref{eqInt.4}) satisfy
condition (\ref{eqTig.1}); for example, the representation $\left(
T_{i}\right)  _{i=0}^{N-1}$ defined from the constant matrix function
$\mathbb{T}\ni z\mapsto I_{N}=%
\begin{pmatrix}
1 &  & \phantom{0}\raisebox{-6pt}[0pt][0pt]{\llap{\Huge $0$}}\\
& \smash{\ddots}\rule{0pt}{10pt} & \\
\raisebox{0pt}[0pt][0pt]{\rlap{\Huge $0$}}\phantom{0} &  & 1
\end{pmatrix}
\in\mathrm{U}_{N}\left(  \mathbb{C}\right)  $ clearly does not satisfy
(\ref{eqTig.1}). Yet as we show, we may use this simple representation as a
base-point for a comparison with all other representations. Specifically, we
have the following lemma.

\begin{lemma}
\label{LemTig.3}If $S_{i}=S_{i}^{\left(  A\right)  }$ is any representation
defined from some matrix function $\mathbb{T}\ni z\mapsto A\left(  z\right)
\in\mathrm{U}_{N}\left(  \mathbb{C}\right)  $, then
\begin{equation}
\left(  T_{j}^{\ast}S_{i}^{\left(  A\right)  }\right)  f\left(  z\right)
=A_{i,j}\left(  z\right)  f\left(  z\right)  ,\qquad i,j=0,\dots,N-1,\;f\in
L^{2}\left(  \mathbb{T}\right)  \mathpunct{;} \label{eqTig.7}%
\end{equation}
i.e., when $i$, $j$ are given, then the operator $T_{j}^{\ast}S_{i}^{\left(
A\right)  }$ is a multiplication operator on $L^{2}\left(  \mathbb{T}\right)
$, in fact multiplication by the matrix entry $A_{i,j}\left(  z\right)  $ of
the unitary matrix $A\left(  z\right)  $.
\end{lemma}

\begin{proof}
Apply formula (\ref{eqInt.4}) and the fact that the operators $T_{j}$ and
$T_{j}^{\ast}$ are given on $L^{2}\left(  \mathbb{T}\right)  $ as
\begin{align}
T_{j}f\left(  z\right)   &  =z^{j}f\left(  z^{N}\right)  \text{,\qquad
and}\label{eqTig.8}\\
\left(  T_{j}^{\ast}f\right)  \left(  z\right)   &  =\frac{1}{N}%
\sum_{\substack{w\in\mathbb{T}\\w^{N}=z}}w^{-j}f\left(  w\right)  ,\qquad f\in
L^{2}\left(  \mathbb{T}\right)  . \label{eqTig.8bis}%
\end{align}

\end{proof}

As a corollary we get the following formula for the adjoint $S_{i}^{\left(
A\right)  \,\ast}$ in general:%
\begin{equation}
S_{i}^{\left(  A\right)  \,\ast}f\left(  z\right)  =\sum_{j=0}^{N-1}%
\overline{A_{i,j}\left(  z\right)  }\,T_{j}^{\ast}f\left(  z\right)  .
\label{eqTig.10}%
\end{equation}

When $A$ is given we introduce the subspace%
\begin{equation}
\mathcal{L}=\mathcal{L}^{\left(  A\right)  }=\bigvee_{I\in\mathcal{I}\left(
N\right)  }\left[  S_{I}^{\left(  A\right)  \,\ast}\openone\right]  ,
\label{eqTig.11}%
\end{equation}
where $\openone=e_{0}$ is the constant function on $\mathbb{T}$, consistent
with the terminology $e_{n}\left(  z\right)  =z^{n}$, $n\in\mathbb{Z}$. If the
matrix entries in (\ref{eqTig.7}) are Fourier polynomials, it is clear that
$\mathcal{L}^{\left(  A\right)  }$ is a finite-dimensional co-invariant
subspace, and therefore also pure. In the next result, we give a necessary and
sufficient condition for $\mathcal{L}^{\left(  A\right)  }$ to be a pure
co-invariant subspace for the wavelet representation $\left(  S_{j}^{\left(
A\right)  }\right)  _{j=0}^{N-1}$ on $L^{2}\left(  \mathbb{T}\right)  $.

\begin{lemma}
\label{LemTig.4}Let $\mathbb{T}\ni z\mapsto A\left(  z\right)  \in
\mathrm{U}_{N}\left(  \mathbb{C}\right)  $ be a measurable unitary
matrix-valued function, and let $\left(  S_{j}\right)  =\left(  S_{j}^{\left(
A\right)  }\right)  $ be the corresponding representation of $\mathcal{O}_{N}$
on the Hilbert space $L^{2}\left(  \mathbb{T}\right)  $. Then the following
three conditions are equivalent. \textup{(}We consider $0\leq j<N$ and
functions in $L^{2}\left(  \mathbb{T}\right)  $.\textup{)}

\begin{enumerate}
\item \label{LemTig.4(1)}%
\begin{equation}
S_{j}f=g. \label{eqTig.12}%
\end{equation}

\item \label{LemTig.4(2)}%
\[
S_{j}^{\ast}g=f\text{\qquad and\qquad}\left\Vert g\right\Vert =\left\Vert
f\right\Vert .
\]

\item \label{LemTig.4(3)}%
\begin{equation}
T_{i}^{\ast}g=f\cdot A_{j,i}\text{\qquad for all }i. \label{eqTig.13}%
\end{equation}

\end{enumerate}

Here each of the identities in \textup{(\ref{LemTig.4(1)})}%
--\textup{(\ref{LemTig.4(3)})} is taken in the pointwise sense, i.e., identity
for the functions on $\mathbb{T}$ pointwise $\mathrm{a.e.}$ with respect to
Haar measure on $\mathbb{T}$. The product on the right-hand side in
\textup{(\ref{eqTig.13})} is $f\left(  z\right)  A_{j,i}\left(  z\right)  $
$\mathrm{a.e.}\;z\in\mathbb{T}$.
\end{lemma}

\begin{proof}
(\ref{LemTig.4(1)})${}\Rightarrow{}$(\ref{LemTig.4(2)}): This is clear since
$S_{j}$ is an isometry. Hence $\left\Vert f\right\Vert =\left\Vert
g\right\Vert $, and $S_{j}^{\ast}g=S_{j}^{\ast}S_{j}f=f$, which is the
combined assertion in (\ref{LemTig.4(2)}).

(\ref{LemTig.4(2)})${}\Rightarrow{}$(\ref{LemTig.4(3)}): Assuming
(\ref{LemTig.4(2)}), and using (\ref{eqTig.10}), we get%
\begin{equation}
\sum_{i=0}^{N-1}\overline{A_{j,i}\left(  z\right)  }\,T_{i}^{\ast}g\left(
z\right)  =f\left(  z\right)  ,\qquad z\in\mathbb{T}. \label{eqTig.14}%
\end{equation}
Using unitarity of the matrix function $A$, and the Schwarz inequality for the
Hilbert space $\mathbb{C}^{N}$, we get the pointwise estimate%
\begin{equation}
\left\vert f\left(  z\right)  \right\vert ^{2}\leq\sum_{i=0}^{N-1}\left\vert
T_{i}^{\ast}g\left(  z\right)  \right\vert ^{2},\qquad\mathrm{a.e.}%
\;z\in\mathbb{T}. \label{eqTig.15}%
\end{equation}
Integration of this over $\mathbb{T}$ with respect to Haar measure yields%
\begin{equation}
\left\Vert f\right\Vert ^{2}\leq\sum_{i=0}^{N-1}\left\Vert T_{i}^{\ast
}g\right\Vert ^{2}=\left\Vert g\right\Vert ^{2}. \label{eqTig.16}%
\end{equation}
But the second condition in (\ref{LemTig.4(2)}) then states that we have
equality in Schwarz's inequality. First we have it in the vector form
(\ref{eqTig.16}). But this means that%
\[
\int_{\mathbb{T}}\left(  \sum_{i=0}^{N-1}\left\vert T_{i}^{\ast}g\left(
z\right)  \right\vert ^{2}-\left\vert f\left(  z\right)  \right\vert
^{2}\right)  \,d\mu\left(  z\right)  =0.
\]
In view of (\ref{eqTig.15}), this means that in fact, (\ref{eqTig.15}) is an
$\mathrm{a.e.}$ identity, i.e., that%
\[
\left\vert f\right\vert ^{2}=\sum_{i=0}^{N-1}\left\vert T_{i}^{\ast
}g\right\vert ^{2}\qquad\mathrm{a.e.}\text{ on }\mathbb{T},
\]
and that therefore%
\[
\left\vert \sum_{i=0}^{N-1}\overline{A_{j,i}\left(  z\right)  }\,T_{i}^{\ast
}g\left(  z\right)  \right\vert ^{2}=\sum_{i=0}^{N-1}\left\vert T_{i}^{\ast
}g\left(  z\right)  \right\vert ^{2}\text{\qquad for }\mathrm{a.e.}%
\;z\in\mathbb{T}.
\]
Hence there is a function $h_{j}$ on $\mathbb{T}$ such that $T_{i}^{\ast
}g=h_{j}A_{j,i}$. But an application of (\ref{LemTig.4(2)}) and $\sum
_{i=0}^{N-1}T_{i}T_{i}^{\ast}=I_{L^{2}\left(  \mathbb{T}\right)  }$ shows that
$h_{j}=f$, which is the desired conclusion (\ref{LemTig.4(3)}), i.e., the
formula (\ref{eqTig.13}).

(\ref{LemTig.4(3)})${}\Rightarrow{}$(\ref{LemTig.4(1)}): If (\ref{LemTig.4(3)}%
) holds, we get%
\begin{align*}
S_{j}f\left(  z\right)   &  =\sum_{i=0}^{N-1}A_{j,i}\left(  z^{N}\right)
T_{i}f\left(  z\right)  =\sum_{i=0}^{N-1}z^{i}T_{i}^{\ast}g\left(
z^{N}\right) \\
&  =\sum_{i=0}^{N-1}T_{i}T_{i}^{\ast}g\left(  z\right)  =g\left(  z\right)  .
\end{align*}

\end{proof}

We note two consequences deriving from the condition (\ref{eqTig.1}). It is a
condition on the given measurable matrix function $\mathbb{T}\ni z\mapsto
A\left(  z\right)  \in\mathrm{U}_{N}\left(  \mathbb{C}\right)  $, and
therefore on the corresponding representation $\left(  S_{j}^{\left(
A\right)  }\right)  _{j=0}^{N-1}$ of $\mathcal{O}_{N}$. This representation
acts on the Hilbert space $L^{2}\left(  \mathbb{T}\right)  $. But the wavelet
system $\varphi$, $\psi_{j}$ from (\ref{eqInt.1}) relates to the line
$\mathbb{R}$, and not directly to $\mathbb{T}$. Indeed, condition
(\ref{eqTig.1}) ensures that the wavelet functions of the system $\varphi$,
$\psi_{j}$, derived from $A$, are in $L^{2}\left(  \mathbb{R}\right)  $, and
(\ref{eqTig.1}) is called the frequency-subband condition. The functions
\begin{equation}
m_{j}^{\left(  A\right)  }\left(  z\right)  =\sum_{k=0}^{N-1}A_{j,k}\left(
z^{N}\right)  z^{k} \label{eqTig.17}%
\end{equation}
are called subband filters: $m_{0}^{\left(  A\right)  }$ is the low-pass
filter, and the others $m_{j}^{\left(  A\right)  }$, $j\geq1$, are the
higher-pass filter bands.

\begin{notations}
\label{NotTigNew.5}Condition (\ref{eqTig.1}) gives the distribution of the $N$
cases with
probabilities $\frac{1}{N}\left\vert m_{j}^{\left(  A\right)  }\left(
\,\cdot\,\right)  \right\vert ^{2}$ on the frequencies $\frac{j}{N}$,
$j=0,1,\dots,N-1$, which represent the bands. Recall if $\rho=e^{i2\pi/N}$,
then $\left\{  \rho^{j}:0\leq j<N\right\}  $ are the $N$'th roots of unity,
i.e., $\left(  \rho^{j}\right)  ^{N}=1$. The frequency passes for the bands
$0,\frac{1}{N},\dots,\frac{N-1}{N}$, referring to the low-pass filter, are
$(\underbrace{1,}_{\text{pass}}\underbrace{0,0,\dots,0}_{\text{halt}})$, and
similarly%
\begin{equation}
\frac{1}{N}\left\vert m_{j}^{\left(  A\right)  }\left(  \rho^{k}\right)
\right\vert ^{2}=\delta_{j,k}. \label{eqTig.18}%
\end{equation}
In fact, these conditions (\ref{eqTig.18}) are equivalent to the single matrix
condition (\ref{eqTig.1}) for $A$.

We also note that (\ref{eqTig.1}) implies that each one of the $N$ isometries
$S_{j}^{\left(  A\right)  }$ on $L^{2}\left(  \mathbb{T}\right)  $ is a shift,
i.e., that%
\begin{equation}
\lim_{n\rightarrow\infty}S_{j}^{\left(  A\right)  \ast\,n}=0. \label{eqTig.19}%
\end{equation}
This conclusion, while nontrivial, is contained in the result Theorem 3.1 in
\cite{BrJo97b}. Note that each of the $N$ shift operators has infinite
multiplicity in the sense of Proposition \ref{ProPur.2} above. Recall if
$S\colon\mathcal{H}\rightarrow\mathcal{H}$ is a shift in a Hilbert space
$\mathcal{H}$, then the multiplicity space is $\mathcal{W}_{S}:=\left(
S\mathcal{H}\right)  ^{\perp}=\ker\left(  S^{\ast}\right)  $: specifically,
$\mathcal{H}\simeq
\smash[b]{\sideset{}{^{\smash{\oplus}}}{\textstyle\sum}\limits_{n=0}%
^{\infty}}\mathcal{W}_{S}$ with $S$ represented as%
\begin{equation}
\left(  x_{0},x_{1},\dots\right)  \overset{\hat{S}}{\longmapsto}\left(
0,x_{0},x_{1},\dots\right)  \label{eqTig.20}%
\end{equation}
where $x_{0},x_{1},\dots\in\mathcal{W}_{S}$. For the particular application to
$S_{j}^{\left(  A\right)  }$, we have%
\begin{equation}
\mathcal{W}_{S_{j}^{\left(  A\right)  }}%
=\sideset{}{^{\oplus}}{\textstyle\sum}\limits_{k\neq j}S_{k}^{\left(
A\right)  }\left(  L^{2}\left(  \mathbb{T}\right)  \right)  . \label{eqTig.21}%
\end{equation}

\end{notations}

The next result is a corollary of Proposition \ref{ProPur.2} and the
discussion above.

\begin{corollary}
\label{CorTig.5}\textup{(A Dichotomy.)} Let the matrix function $\mathbb{T}\ni
z\mapsto A\left(  z\right)  \in\mathrm{U}_{N}\left(  \mathbb{C}\right)  $,
$N>1$, satisfy \textup{(\ref{eqTig.1})}, and let $\left(  S_{j}^{\left(
A\right)  }\right)  _{j=0}^{N-1}$ be the corresponding representation of
$\mathcal{O}_{N}$ on $L^{2}\left(  \mathbb{T}\right)  $. If $\mathcal{L}%
\subset L^{2}\left(  \mathbb{T}\right)  $ is a co-invariant subspace, then the
following two conditions are equivalent.

\begin{enumerate}
\item \label{CorTig.5(1)}$\mathcal{L}$ is pure.

\item \label{CorTig.5(2)}$\mathcal{L}\neq L^{2}\left(  \mathbb{T}\right)  $.
\end{enumerate}
\end{corollary}

\begin{proof}
Clearly (\ref{CorTig.5(1)})${}\Rightarrow{}$(\ref{CorTig.5(2)}), and it is
immediate from Definition \ref{DefPur.1} that every finite-dimensional
co-invariant subspace $\mathcal{L}$ is not pure. That follows since, if
$f\in\mathcal{L}$ satisfies $S^{\infty}f\subset\mathcal{L}$, then for each
$k\in\mathbb{Z}_{+}$, the family $\left\{  S_{I}f:I\in\mathcal{I}\left(
N\right)  ,\;\left\vert I\right\vert =k\right\}  \subset S^{\infty}%
f\subset\mathcal{L}$ consists of orthogonal vectors. If $f\neq0$, then there
are $N^{k}$ such vectors. Hence the result holds whenever the entries in $A$
are Fourier polynomials; see also \cite{JoKr02}. We now turn to
(\ref{CorTig.5(2)})${}\Rightarrow{}$(\ref{CorTig.5(1)}) in the general case:
Let $A$ satisfy (\ref{eqTig.1}), and let $\mathcal{L}$ be a co-invariant
subspace, referring to the representation of $\mathcal{O}_{N}$ given by
$\left(  S_{j}^{\left(  A\right)  }\right)  _{j=0}^{N-1}$. Then $\mathcal{L}%
^{\perp}=L^{2}\left(  \mathbb{T}\right)  \ominus\mathcal{L}$ is invariant for
each of the $N$ shift operators $S_{j}^{\left(  A\right)  }$, $j=0,\dots,N-1$.
Now suppose $\mathcal{L}\neq L^{2}\left(  \mathbb{T}\right)  $, or
equivalently that $\mathcal{L}^{\perp}\neq0$. Setting $j=0$, we see that there
is a unitary operator-valued function $\mathbb{T}\ni z\mapsto u_{0}\left(
z\right)  \in\mathrm{U}\left(  \mathcal{W}_{S_{0}^{\left(  A\right)  }%
}\right)  $ such that
\begin{equation}
\mathcal{L}^{\perp}=\left\{  u_{0}h_{+}:h_{+}\in H_{+}\left(  \mathcal{W}%
_{S_{0}^{\left(  A\right)  }}\right)  \right\}  \label{eqTig.22}%
\end{equation}
where $\mathcal{W}_{S_{0}^{\left(  A\right)  }}$ is given by (\ref{eqTig.21}).

If $f\in\mathcal{L}$ satisfies $S^{\left(  A\right)  \,\infty}f\subset
\mathcal{L}$, then $f,S_{0}^{\left(  A\right)  }f,S_{0}^{\left(  A\right)
\,2}f,\dots,S_{0}^{\left(  A\right)  \,n}f,\dots\in\mathcal{L}$. Relative to
the representation (\ref{eqTig.20}), applied to $S_{0}^{\left(  A\right)  }$,
we get%
\[
\hat{f},z\hat{f},z^{2}\hat{f},\dots,z^{n}\hat{f},\dots\in\mathcal{L}.
\]
But then we get $f=0$ by an application of Proposition \ref{ProPur.2} to
$S_{0}^{\left(  A\right)  }$. Hence, $\mathcal{L}$ is pure; see Definition
\ref{DefPur.1}. Note that \cite[Theorem 3.1]{BrJo97b} was used as well. This
result implies that $S_{0}^{\left(  A\right)  }$ is a shift, and so it has the
representation (\ref{eqTig.20}).
\end{proof}

We now turn to the space of Lipschitz functions on $\mathbb{T}$. Via the
coordinate $z=e^{-i\xi}$, $\xi\in\mathbb{R}$, we identify functions on
$\mathbb{T}$ with $2\pi$-periodic functions on $\mathbb{R}$ and we define the
\emph{Lipschitz space} $\operatorname*{Lip}\nolimits_{1}$ by%
\begin{equation}
\left\Vert f\right\Vert _{\operatorname*{Lip}\nolimits_{1}}:=\left\vert
f\left(  0\right)  \right\vert +\sup_{-\pi\leq\xi<\eta<\pi}\frac{\left\vert
f\left(  \xi\right)  -f\left(  \eta\right)  \right\vert }{\left\vert \xi
-\eta\right\vert }<\infty. \label{eqTig.23}%
\end{equation}
A matrix function is said to be Lipschitz if its matrix entries are in
$\operatorname*{Lip}\nolimits_{1}$.

\begin{theorem}
\label{ThmTig.6}A matrix function, $\mathbb{T}\ni z\mapsto A\left(  z\right)
\in\mathrm{U}_{N}\left(  \mathbb{C}\right)  $ is given. We assume it is in the
Lipschitz class, and that it satisfies \textup{(\ref{eqTig.1})}. Consider the
co-invariant subspace $\mathcal{L}^{\left(  A\right)  }\subset L^{2}\left(
\mathbb{T}\right)  $ defined from the corresponding representation
$\smash[b]{\left(
S_{j}^{\left(  A\right)  }\right)  _{j=0}^{N-1}}$ of $\mathcal{O}_{N}$ as
follows:%
\begin{equation}
\mathcal{L}^{\left(  A\right)  }:=\bigvee\left[  S_{I}^{\left(  A\right)
\,\ast}e_{0}:I\in\mathcal{I}\left(  N\right)  \right]  \label{eqTig.24}%
\end{equation}
where $e_{0}=\openone$ is the constant function $1$ on $\mathbb{T}$. Then
$\mathcal{L}^{\left(  A\right)  }$ is pure.
\end{theorem}

\begin{proof}
In view of Corollary \ref{CorTig.5}, it is enough to show that $\mathcal{L}%
^{\left(  A\right)  }\neq L^{2}\left(  \mathbb{T}\right)  $, or equivalently
that $\left(  \mathcal{L}^{\left(  A\right)  }\right)  ^{\perp}\neq0$. The
argument is indirect. If $\mathcal{L}^{\left(  A\right)  }=L^{2}\left(
\mathbb{T}\right)  $, we get a contradiction as follows (we have suppressed
the superscript $A$ in the notation): The set of vectors $\left\{  S_{I}%
^{\ast}e_{0}:I\in\mathcal{I}\left(  N\right)  \right\}  $ is relatively
compact in $C\left(  \mathbb{T}\right)  $ by Arzel\`{a}-Ascoli. To see this,
we use the Lipschitz property, and the formula%
\[
S_{I}^{\ast}e_{0}\left(  z\right)  =\frac{1}{N^{k}}\sum_{w^{N^{k}}=z}%
\overline{m_{I}^{\left(  k\right)  }\left(  w\right)  }%
\]
where $m_{I}^{\left(  k\right)  }\left(  z\right)  =m_{i_{1}}\left(  z\right)
m_{i_{2}}\left(  z^{N}\right)  \cdots m_{i_{k}}\left(  z^{N^{k-1}}\right)  $,
$I=\left(  i_{1},\dots,i_{k}\right)  $. Further note that $\left(  \sqrt
{N}S_{0}^{\ast}f\right)  \left(  1\right)  =f\left(  1\right)  $ for all
Lipschitz functions $f$, and that the analogous conditions hold for
$S_{1}^{\ast},\dots,S_{N-1}^{\ast}$. Since $m_{j}\left(  \rho^{k}\right)
=\delta_{j,k}\sqrt{N}$, there is a Lipschitz solution $f$ to the following
system of equations:
\begin{equation}
f\left(  1\right)  =1,\qquad S_{0}^{\ast}f=\frac{1}{\sqrt{N}}f,\qquad
S_{j}^{\ast}f=0,\;j\geq1. \label{eqTig.25}%
\end{equation}
Using $f=\sum_{j}S_{j}S_{j}^{\ast}f=\frac{1}{\sqrt{N}}S_{0}f$, it follows that
$S_{0}$ has $\sqrt{N}$ as eigenvalue, contradicting that $S_{0}$ is isometric.
The contradiction proves that $\mathcal{L}^{\left(  A\right)  \,\perp}\neq0$.
\end{proof}

In the next section, we give additional details on the existence question for
the eigenvalue problems related to the operators $S_{j}^{\left(  A\right)
\,\ast}$ in the case when the matrix function $A$ is assumed to be Lipschitz.

\section{\label{Fin}Finite dimensions}

In this section we offer a construction of a finite-dimensional co-invariant
(nonzero, and nontrivial) subspace for the representation of $\mathcal{O}_{N}$
on $L^{2}\left(  \mathbb{T}\right)  $ which is associated with a
multiresolution wavelet of scale $N$. It is both natural and optimal with
respect to the conditions of Sections \ref{Non}--\ref{Tig}. The setting is as
in the previous section: Recall, a Lipschitz mapping $\mathbb{T}\ni z\mapsto
A\left(  z\right)  \in\mathrm{U}_{N}\left(  \mathbb{C}\right)  $ is given, and
it is assumed that the subbands are ordered according to (\ref{eqTig.1}),
i.e., that $A\left(  1\right)  $ is the Hadamard matrix, or equivalently that%
\begin{equation}
m_{j}\left(  \rho_{N}^{k}\right)  =\delta_{j,k}\sqrt{N},\qquad j,k=0,\dots
,N-1, \label{eqFin.1}%
\end{equation}
where $\rho_{N}:=\exp\left(  i\frac{2\pi}{N}\right)  $, and
\begin{equation}
m_{j}\left(  z\right)  =\sum_{k=0}^{N-1}A_{j,k}\left(  z^{N}\right)  z^{k}.
\label{eqFin.2}%
\end{equation}
Then the operators $S_{j}=S_{j}^{\left(  A\right)  }$ on $L^{2}\left(
\mathbb{T}\right)  $ are%
\begin{equation}
\left(  S_{j}f\right)  \left(  z\right)  =m_{j}\left(  z\right)  f\left(
z^{N}\right)  ,\qquad j=0,\dots,N-1,\;f\in L^{2}\left(  \mathbb{T}\right)
,\;z\in\mathbb{T}. \label{eqFin.3}%
\end{equation}

\begin{theorem}
\label{ThmFin.1}Let $\left(  S_{j}\right)  _{j=0}^{N-1}$ be a representation
as in \textup{(\ref{eqFin.3})} determined by a Lipschitz system and subject to
conditions \textup{(\ref{eqFin.1})} and \textup{(\ref{eqFin.2})}. Then the
following two conditions are equivalent.

\begin{enumerate}
\item \label{ThmFin.1(1)}There is a finite-dimensional co-invariant subspace
$\mathcal{L}\subset L^{2}\left(  \mathbb{T}\right)  $ which contains the
solutions $f$ to the following affine conditions:%
\begin{equation}
f\in\operatorname*{Lip}\nolimits_{1},\qquad f\left(  1\right)  =1,\qquad
\sqrt{N}S_{0}^{\ast}f=f. \label{eqFin.4}%
\end{equation}

\item \label{ThmFin.1(2)}There is a finite constant $K$ such that, for all
$f\in\operatorname*{Lip}\nolimits_{1}$,%
\[
\sup_{J\in\mathcal{I}\left(  N\right)  }\left\vert N^{\frac{\left\vert
J\right\vert }{2}}\left(  S_{J}^{\ast}f\right)  \left(  1\right)  \right\vert
\leq K\left\Vert f\right\Vert .
\]

\end{enumerate}

\noindent The affine dimension of the convex set \textup{(\ref{eqFin.4})} of
Lipschitz functions is at least one.
\end{theorem}

\begin{proof}
We begin with the conditions (\ref{eqFin.4}). We will show that there is a
well defined linear operator $T$ on $\operatorname{Lip}_{1}$ with
finite-dimensional range (dimension at least one) such that
\begin{equation}
\left(  Tf\right)  \left(  1\right)  =f\left(  1\right)  \text{\qquad for all
}f\in\operatorname*{Lip}\nolimits_{1} \label{eqFin.5}%
\end{equation}
and%
\begin{equation}
\left(  Tf\right)  \left(  z\right)  =\lim_{k\rightarrow\infty}\left(
\sqrt{N}S_{0}^{\ast}\right)  ^{k}f\left(  z\right)  ,\qquad f\in
\operatorname*{Lip}\nolimits_{1}, \label{eqFin.6}%
\end{equation}
where the limit in (\ref{eqFin.6}) is uniform for $z\in\mathbb{T}$.

\begin{lemma}
\label{LemFin.2}Let the operators $S_{j}$, $j=0,\dots,N-1$, be as specified
above in \textup{(\ref{eqFin.3})}, i.e., the Lipschitz property is assumed, as
is \textup{(\ref{eqFin.1})}. For bounded functions $f$ on $\mathbb{T}$, set
$\left\Vert f\right\Vert :=\sup_{z\in\mathbb{T}}\left\vert f\left(  z\right)
\right\vert $; and if $f$ is differentiable, set
\begin{equation}
\tilde{f}\left(  x\right)  :=f\left(  e^{-i2\pi x}\right)  , \label{eqFin.7}%
\end{equation}
and $f^{\prime}:=\tilde{f}^{\prime}$. Finally, let%
\begin{equation}
M_{1}:=N^{-\frac{1}{2}}\max_{0\leq j<N}\left\Vert m_{j}^{\prime}\right\Vert .
\label{eqFin.8}%
\end{equation}
Then we have the following estimate:%
\begin{multline}
N^{\frac{k}{2}}\left\Vert \left(  S_{j_{k}}^{\ast}\cdots S_{j_{2}}^{\ast
}S_{j_{1}}^{\ast}f\right)  ^{\prime}\right\Vert \leq N^{-\frac{k}{2}%
}\left\Vert f^{\prime}\right\Vert +M_{1}\left\Vert f\right\Vert
\label{eqFin.9}\\
\text{for all }k\in\mathbb{Z}_{+}\text{, all }f\in\operatorname*{Lip}%
\nolimits_{1}\text{,}\\
\text{and all multi-indices }J=\left(  j_{1},j_{2},\dots,j_{k}\right)
\in\mathcal{I}_{k}\left(  N\right)  .
\end{multline}

\end{lemma}

\begin{proof}
Using (\ref{eqFin.7}), we shall pass freely between any of the four equivalent
formulations, functions on $\mathbb{T}$, functions on $\left[  0,1\right)  $,
functions on $\mathbb{R}\diagup\mathbb{Z}$, or one-periodic functions on
$\mathbb{R}$, omitting the distinction between $f$ and $\tilde{f}$ in
(\ref{eqFin.7}). With the multi-index notation in (\ref{eqFin.9}), we set%
\begin{equation}
m_{J}\left(  z\right)  :=m_{j_{1}}\left(  z\right)  m_{j_{2}}\left(
z^{N}\right)  \cdots m_{j_{k}}\left(  z^{N^{k-1}}\right)  ,\qquad
z\in\mathbb{T}. \label{eqFin.10}%
\end{equation}
Setting $S_{J}:=S_{j_{1}}\cdots S_{j_{k}}$, and $S_{J}^{\ast}:=S_{j_{k}}%
^{\ast}\cdots S_{j_{2}}^{\ast}S_{j_{1}}^{\ast}$, we find%
\begin{align}
\left(  S_{J}f\right)  \left(  z\right)  &= m_{J}\left(  z\right)  f\left(
z^{N^{k}}\right)  , \label{eqFin.11a}\\
\intertext{and}\left(  S_{J}^{\ast}f\right)  \left(
z\right)  &= \frac{1}{N^{k}}\!\!\sum_{\substack{w\in\mathbb{T}\\w^{N^{k}}%
=z}}\!\!\overline{m_{J}\left(  w\right)  }\,f\left(  w\right)  .
\label{eqFin.11b}%
\end{align}
The sum in (\ref{eqFin.11b}) contains $N^{k}$ terms, which is evident from the
rewrite in the form below, using instead one-periodic functions on
$\mathbb{R}$, $x\in\mathbb{R}$:%
\begin{equation}
\left(  S_{J}^{\ast}f\right)  \left(  x\right)  =\frac{1}{N^{k}}%
\!\!\sum_{\substack{y\in\mathbb{R}\\N^{k}y\equiv x\bmod1}}\!\!\overline
{m_{J}\left(  y\right)  }\,f\left(  y\right)  . \label{eqFin.12}%
\end{equation}
In this form $m_{J}\left(  y\right)  =m_{j_{1}}\left(  y\right)  m_{j_{2}%
}\left(  Ny\right)  \cdots m_{j_{k}}\left(  N^{k-1}y\right)  $ (compare with
(\ref{eqFin.10})), and the points $y$ may be represented as%
\begin{equation}
y=\frac{x+l_{0}+l_{1}N+l_{2}N^{2}+\dots+l_{k-1}N^{k-1}}{N^{k}},
\label{eqFin.13}%
\end{equation}
with the integers $l_{0}$, $l_{1}$, $\dots$, $l_{k-1}$ taking values over the
$\operatorname*{mod}N$ residue classes $0,1,\dots,N-1$. Given this, it is
clear how the general form of (\ref{eqFin.9}) follows from the case $k=1$. We
will do $k=1$, and leave the induction and the multi-index gymnastics to the
reader. Using (\ref{eqFin.12})--(\ref{eqFin.13}) we have%
\begin{multline}
N^{\frac{1}{2}}\left(  S_{j}^{\ast}f\right)  ^{\prime}\left(  x\right)
=\frac{1}{N\sqrt{N}}\sum_{l=0}^{N-1}\bar{m}_{j}\left(  \frac{x+l}{N}\right)
f^{\prime}\left(  \frac{x+l}{N}\right) \label{eqFin.14}\\
+\frac{1}{N\sqrt{N}}\sum_{l=0}^{N-1}\bar{m}_{j}^{\prime}\left(  \frac{x+l}%
{N}\right)  f\left(  \frac{x+l}{N}\right)  .
\end{multline}
For the individual terms on the right-hand side, we use Schwarz's inequality,
as follows: First,%
\begin{align*}
&  \frac{1}{N}\left\vert \sum_{l=0}^{N-1}\bar{m}_{j}\left(  \frac{x+l}%
{N}\right)  f^{\prime}\left(  \frac{x+l}{N}\right)  \right\vert \\
&  \qquad\qquad\leq
\vphantom{\underbrace{\frac{1}{N}\sum_{l=0}^{N-1}\left\vert m_{j}\left( \frac{x+l}{N}\right) \right\vert ^{2}}_{=1}}\left(
\smash{\underbrace{\frac{1}{N}\sum_{l=0}^{N-1}\left\vert m_{j}\left( \frac{x+l}{N}\right) \right\vert ^{2}}_{=1}}\vphantom{\frac{1}{N}\sum_{l=0}^{N-1}\left\vert m_{j}\left( \frac{x+l}{N}\right) \right\vert ^{2}}\right)
^{\frac{1}{2}}\left(  \frac{1}{N}\sum_{l=0}^{N-1}\left\vert f^{\prime}\left(
\frac{x+l}{N}\right)  \right\vert ^{2}\right)  ^{\frac{1}{2}}\\
&  \qquad\qquad\leq\left\Vert f^{\prime}\right\Vert ;
\end{align*}
and second,%
\begin{align*}
&  \frac{1}{N}\left\vert \sum_{l=0}^{N-1}\bar{m}_{j}^{\prime}\left(
\frac{x+l}{N}\right)  f\left(  \frac{x+l}{N}\right)  \right\vert \\
&  \qquad\qquad\leq\left(  \frac{1}{N}\sum_{l=0}^{N-1}\left\vert m_{j}%
^{\prime}\left(  \frac{x+l}{N}\right)  \right\vert ^{2}\right)  ^{\frac{1}{2}%
}\left(  \frac{1}{N}\sum_{l=0}^{N-1}\left\vert f\left(  \frac{x+l}{N}\right)
\right\vert ^{2}\right)  ^{\frac{1}{2}}\\
&  \qquad\qquad\leq\sqrt{N}M_{1}\left\Vert f\right\Vert ,
\end{align*}
where the terms on the right are given in (\ref{eqFin.8}) and the discussion
in that paragraph.

Introducing the last two estimates back into (\ref{eqFin.14}), we get%
\[
N^{\frac{1}{2}}\left\vert \left(  S_{j}^{\ast}f\right)  ^{\prime}\left(
x\right)  \right\vert \leq N^{-\frac{1}{2}}\left\Vert f^{\prime}\right\Vert
+M_{1}\left\Vert f\right\Vert ,
\]
which is the desired estimate.
\end{proof}

Note that if the functions $f$ in (\ref{eqFin.9}) are restricted by
$\left\Vert f\right\Vert \leq1$, then all the functions $N^{\frac{k}{2}}%
S_{J}^{\ast}f$ are contained in a compact subset in the Banach space $C\left(
\mathbb{T}\right)  $. This follows from \cite[vol.~I, p.~245, Theorem
IV.3.5]{DuSc}; see also \cite[Section 4]{IoMa50}, or \cite{Bal00}. The
conclusion is that there is a fixed finite-dimensional subspace $\mathcal{L}$
which contains all the functions obtained as limits of the terms $N^{\frac
{k}{2}}S_{J}^{\ast}f$ as the multi-indices $J=\left(  j_{1},\dots
,j_{k}\right)  $ vary. Introducing $\mathcal{I}_{\infty}\left(  N\right)
=\left\{  0,1,\dots,N-1\right\}  ^{\mathbb{Z}_{+}}$, for each $\xi=\left(
j_{1},j_{2},\dots\right)  \in\mathcal{I}_{\infty}\left(  N\right)  $ we may
pass to infinity via an ultrafilter $u\left(  \xi\right)  $, i.e.,%
\begin{equation}
T_{\xi}f:=\lim_{\substack{u\left(  \xi\right)  \\\xi=\left(  J_{\rightarrow
}\ast\cdots\ast\cdots\right)  }}N^{\frac{\left\vert J\right\vert }{2}}%
S_{J}^{\ast}f, \label{eqFin.15}%
\end{equation}
and we note that all the operators $T_{\xi}$, $\xi\in\mathcal{I}_{\infty
}\left(  N\right)  $, have their range contained in $\mathcal{L}$. If
$j\in\left\{  0,1,\dots,N-1\right\}  $, then the extension $\left(
j\xi\right)  $ is a point in $\mathcal{I}_{\infty}\left(  N\right)  $ for
every $\xi\in\mathcal{I}_{\infty}\left(  N\right)  $, and%
\begin{equation}
N^{\frac{1}{2}}S_{j}^{\ast}T_{\xi}f=T_{\left(  j\xi\right)  }f\in\mathcal{L}.
\label{eqFin.16}%
\end{equation}
This proves that $\mathcal{L}$ is co-invariant.

For a reference to ultrafilters and compactifications, we suggest
\cite{Encyclo} and \cite{Bourbaki}.

The limit (\ref{eqFin.6}) is covered by this discussion, since the point
$\xi=\left(  0,0,\dots\right)  \in\mathcal{I}_{\infty}\left(  N\right)  $,
so the corresponding operator $T$ in (\ref{eqFin.6}) is just $T_{\left(
0,0,\dots\right)  }$. Hence (\ref{eqFin.6}) is a special case of
(\ref{eqFin.15}). Finally, (\ref{eqFin.5}) follows from%
\begin{equation}
\sqrt{N}\left(  S_{0}^{\ast}f\right)  \left(  0\right)  =\frac{1}{\sqrt{N}%
}\sum_{l=0}^{N-1}\overline{m_{0}\left(  \frac{l}{N}\right)  }\,f\left(
\frac{l}{N}\right)  =f\left(  0\right)  , \label{eqFin.17}%
\end{equation}
which is based on (\ref{eqFin.1}). In the additive form, (\ref{eqFin.1})
implies $m_{0}\left(  0\right)  =\sqrt{N}$, and $m_{0}\left(  \frac{1}%
{N}\right)  =\dots=m_{0}\left(  \frac{N-1}{N}\right)  =0$. Note $z=\exp\left(
-i2\pi\left(  x=0\right)  \right)  =1$ is used in (\ref{eqFin.5}).

Let $f\in\operatorname*{Lip}\nolimits_{1}$, and set $f_{J}=N^{\frac{\left\vert
J\right\vert }{2}}S_{J}^{\ast}f$. We noted that each $f_{J}$ is in
$\operatorname*{Lip}\nolimits_{1}$, and we gave a uniform estimate on the
derivatives, i.e., on $f_{J}^{\prime}\left(  x\right)  $. But the sequence is
also bounded in $C\left(  \mathbb{T}\right)  $, relative to the usual norm
$\left\Vert \,\cdot\,\right\Vert $ on $C\left(  \mathbb{T}\right)  $.
Moreover,
\begin{equation}
\left\Vert f_{J}\right\Vert \leq N^{-\frac{\left\vert J\right\vert }{2}%
}\left\Vert f^{\prime}\right\Vert +\left(  M_{1}+K\right)  \left\Vert
f\right\Vert . \label{eqFin.star1}%
\end{equation}

To see this, first note the estimate from (\ref{ThmFin.1(2)}), i.e.,%
\begin{equation}
\sup_{J\in\mathcal{I}\left(  N\right)  }\left\vert f_{J}\left(  1\right)
\right\vert \leq K\left\Vert f\right\Vert . \label{eqFin.star2}%
\end{equation}
The details are understood best in the additive formulation, i.e., with
$\mathbb{T}\cong\left[  0,1\right)  $ and the identification $z\cong x$ via
$z=\exp\left(  -i2\pi x\right)  $. The first iteration step, starting with
$f$, yields ($x=0$)%
\[
N^{\frac{1}{2}}\left(  S_{0}^{\ast}f\right)  \left(  0\right)  =N^{-\frac
{1}{2}}\sum_{l=0}^{N-1}\bar{m}_{0}\left(  \frac{l}{N}\right)  f\left(
\frac{l}{N}\right)  =f\left(  0\right)  ,
\]
and subsequent steps yield%
\[
N^{\frac{1}{2}}\left(  S_{j}^{\ast}f\right)  \left(  0\right)  =N^{-\frac
{1}{2}}\sum_{l=0}^{N-1}\bar{m}_{j}\left(  \frac{l}{N}\right)  f\left(
\frac{l}{N}\right)  =f\left(  \frac{j}{N}\right)  ,
\]
where conditions (\ref{eqFin.1}) are used. So if $\left\vert J\right\vert =1$,
(\ref{eqFin.star2}) holds, and $\leq$ is equality. Note that in general, if
$\left\vert J\right\vert =k>1$, then
\begin{equation}
F_{J}\left(  0\right)  =N^{-\frac{k}{2}}\sum_{N^{k}x\equiv0\bmod1}%
\overline{m_{J}\left(  x\right)  }\,f\left(  x\right)  , \label{eqFin.star4}%
\end{equation}
where $x$ is an $N$-adic fraction $0\leq x<1$. Note that $x$ varies over the
solutions $N^{k}x\equiv0\bmod1$, and $J$ is fixed. For real values of $x$,
$0\leq x<1$, we have%
\[
\left\vert f_{J}\left(  x\right)  -f_{J}\left(  0\right)  \right\vert
\leqq\int_{0}^{x}\left\vert f_{J}^{\prime}\left(  y\right)  \right\vert
\,dy\leq x\left(  N^{-\frac{k}{2}}\left\Vert f^{\prime}\right\Vert
+M_{1}\left\Vert f\right\Vert \right)  ;
\]
and, using (\ref{eqFin.star2}), we get
\begin{equation}
\left\vert f_{J}\left(  x\right)  \right\vert \leq N^{-\frac{k}{2}}\left\Vert
f^{\prime}\right\Vert +\left(  M_{1}+K\right)  \left\Vert f\right\Vert ,
\label{eqFin.star5}%
\end{equation}
which is the desired estimate (\ref{eqFin.star1}).

If $\xi\in\mathcal{I}_{\infty}\left(  N\right)  $, the operator $T_{\xi}$ is
defined on $f$ as in (\ref{eqFin.15}) by a limit over an ultrafilter,%
\begin{equation}
T_{\xi}f=\!\!\lim_{\substack{u\left(  \xi\right)  \\(\underbrace
{J_{\rightarrow}\ast\cdots\ast\cdots}_{\xi})}}\!\!f_{J}, \label{eqFin.star6}%
\end{equation}
as noted in (\ref{eqFin.15}). Using (\ref{eqFin.star1}), and passing to the
limit, we get $\left\Vert T_{\xi}f\right\Vert \leq\left(  M_{1}+K\right)
\left\Vert f\right\Vert $. Similarly, an application of (\ref{eqFin.star2})
yields $\left\vert \left(  T_{\xi}f\right)  \left(  1\right)  \right\vert \leq
K\left\Vert f\right\Vert $. The fact that the range of $T_{\xi}$ is
finite-dimensional results as noted from Arzel\`{a}-Ascoli in view of the
uniform estimate on the derivatives $f_{J}^{\prime}$.
\end{proof}

\begin{remark}
\label{RemFin.3}Even for the case of the Haar wavelet, both the standard one
where the polyphase matrix $A\left(  z\right)  $ is $A\left(  z\right)  =
\smash[b]{
\begin{pmatrix}
1 & 1\\
1 & -1
\end{pmatrix}
}$, or the stretched Haar wavelet where $A\left(  z\right)  =%
\begin{pmatrix}
1 & z\\
1 & -z
\end{pmatrix}
^{\mathstrut}_{\mathstrut}$,
the space $\left\{  f\in\operatorname*{Lip}\nolimits_{1}:\sqrt{2}%
S_{0}^{\ast}f=f\right\}  $ is of dimension more than $1$. In the first case,
the dimension can be checked to be $2$, and in the second case, it is $3$. The
respective subspaces $\mathcal{L}_{1}$ and $\mathcal{L}_{3}$ are
$\mathcal{L}_{1}=\left[  e_{0},e_{-1}\right]  $ and $\mathcal{L}_{3}=\left[
e_{0},e_{-3},e_{-1}+e_{-2}\right]  $, where $e_{n}\left(  z\right)  =z^{n}$,
$n\in\mathbb{Z}$, $z\in\mathbb{T}$. Two independent functions $f$ in
$\mathcal{L}_{1}$ satisfying $f\left(  1\right)  =1$ are $e_{0}$ and $\frac
{1}{2}\left(  e_{0}+e_{-1}\right)  $. For $\mathcal{L}_{3}$, three such
independent solutions to the system $\sqrt{2}S_{0}^{\ast}f=f$, $f\left(
1\right)  =1$ may be taken to be $e_{0}$, $\frac{1}{2}\left(  e_{0}%
+e_{3}\right)  $, and $\frac{1}{3}\left(  e_{0}+e_{-1}+e_{-2}\right)  $.
\end{remark}

Despite the fact that $\left\{  f\in\operatorname*{Lip}\nolimits_{1}%
:S_{0}^{\ast}f=\left(  1/\sqrt{N}\right)  f\right\}  $ is automatically
finite-dimensional for the representations of Theorem
\ref{ThmFin.1}, we show that for the same class of representations, and for
every $\lambda\in\mathbb{C}$, $\left\vert \lambda\right\vert <1$, the space
\begin{equation}
\left\{  f\in L^{2}\left(  \mathbb{T}\right)  :S_{0}^{\ast}f=\bar{\lambda
}f\right\}  \label{eqFin.18}%
\end{equation}
is infinite-dimensional. If $N=2$, then in fact it is isomorphic to
$\mathcal{W}=S_{1}L^{2}\left(  \mathbb{T}\right)  $. Introducing the
projection $P_{\mathcal{W}}:=S_{1}S_{1}^{\ast}$, it can easily be checked that
the operator $\left(  I-\bar{\lambda}S_{0}\right)  ^{-1}P_{\mathcal{W}}$ is
well defined and maps into $\ker\left(  \bar{\lambda}I-S_{0}^{\ast}\right)  $.
To see that inner products are preserved, note that%
\begin{equation}
\ip{\left( I-\bar{\lambda}S_{0}\right) ^{-1}w_{1}}{\left( I-\bar{\lambda}S_{0}\right) ^{-1}w_{2}}=\frac
{1}{1-\left\vert \lambda\right\vert ^{2}}\ip{w_{1}}{w_{2}} \label{eqFin.19}%
\end{equation}
holds for all $w_{1},w_{2}\in\mathcal{W}$. The easiest way to see this is via
the following known theorem.

\begin{theorem}
\label{ThmFin.5}For every representation $\left(  S_{j}\right)  _{j=0}^{1}$ of
$\mathcal{O}_{2}$ on a Hilbert space $\mathcal{H}$ there is a unique unitary
isomorphism $T\colon H_{+}\left(  \mathcal{W}\right)  \rightarrow\mathcal{H}$
such that $Tw=w$ for $w\in\mathcal{W}$, and%
\begin{equation}
S_{0}T=TM_{z} \label{eqFin.20}%
\end{equation}
if and only if $S_{0}^{\ast\,n}x\underset{n\rightarrow\infty}{\longrightarrow
}0$, $x\in\mathcal{H}$.
\end{theorem}

\begin{proof}
This result is essentially contained in \cite[Theorem 9.1]{BrJo97b}. We will
just give the formula for $T$ and its adjoint $T^{\ast}$, and then leave the
remaining verifications to the reader. $T\colon H_{+}\left(  \mathcal{W}%
\right)  \rightarrow\mathcal{H}$ is
\begin{equation}
T\left(  \sum_{n=0}^{\infty}z^{n}w_{n}\right)  =\sum_{n=0}^{\infty}S_{0}%
^{n}w_{n} \label{eqFinNew.26}%
\end{equation}
where $w_{n}\in\mathcal{W}$ and%
\begin{equation}
\left\Vert \sum_{n=0}^{\infty}z^{n}w_{n}\right\Vert _{H_{+}\left(
\mathcal{W}\right)  }^{2}=\sum_{n=0}^{\infty}\left\Vert w_{n}\right\Vert
_{\mathcal{W}}^{2}, \label{eqFinNew.27}%
\end{equation}
and if $x\in\mathcal{H}$,
\begin{equation}
T^{\ast}x=\sum_{n=0}^{\infty}z^{n}P_{\mathcal{W}}S_{0}^{\ast\,n}x\in
H_{+}\left(  \mathcal{W}\right)  . \label{eqFin.22}%
\end{equation}
The fact that $\left\Vert T^{\ast}x\right\Vert =\left\Vert x\right\Vert $ is
based on the assumption that $S_{0}^{\ast\,n}x\rightarrow0$, i.e., that
$S_{0}$ is a shift.
\end{proof}

Introducing the Szeg\"{o} kernel
\begin{equation}
C_{\lambda}\left(  z\right)  :=\frac{1}{1-\bar{\lambda}z} \label{eqFinNew.29}%
\end{equation}
and the inner function
\begin{equation}
u_{\lambda}\left(  z\right)  =\frac{\lambda-z}{1-\bar{\lambda}z},
\label{eqFinNew.30}%
\end{equation}
we arrive at
\begin{equation}
T\left(  C_{\lambda}\otimes w\right)  =\left(  I-\bar{\lambda}S_{0}\right)
^{-1}w, \label{eqFin.23}%
\end{equation}
which makes (\ref{eqFin.19}) immediate. For the orthocomplement $\left\{
x\in\mathcal{H}:S_{0}^{\ast}x=\bar{\lambda}x\right\}  ^{\perp}$, the
$H_{+}\left(  \mathcal{W}\right)  $ representation is%
\begin{equation}
\left\{  u_{\lambda}\left(  z\right)  F_{+}\left(  z\right)  :F_{+}\in
H_{+}\left(  \mathcal{W}\right)  \right\}  =\left\{  F\in H_{+}\left(
\mathcal{W}\right)  :F\left(  \lambda\right)  =0\right\}  . \label{eqFin.24}%
\end{equation}

\begin{acknowledgements}
We gratefully acknowledge discussions
with David Kribs, Joe Ball and Victor Vinnikov
about representations of the Cuntz
relations, at the 2002 IWOTA
conference at Virginia Tech in
Blacksburg, VA,
and an e-mail with corrections from Ken Davidson.
We thank Brian
Treadway for beautiful typesetting,
for corrections, and for helpful suggestions.
\end{acknowledgements}


\ifx\undefined\bysame
\newcommand{\bysame}{\leavevmode\hbox to3em{\hrulefill}\,}
\fi

\end{document}